\documentclass[12pt]{article}

\usepackage{amsfonts,amsmath}

\def \beq {\begin{eqnarray}}
\def \eeq {\end{eqnarray}}
\def \beqn {\begin{eqnarray*}}
\def \eeqn {\end{eqnarray*}}

\newcommand{\halmos}{\rule{1ex}{1.4ex}}

\newcounter{for}[section]

\newtheorem{itlemma}{Lemma}[section]
\newtheorem{itproposition}[itlemma]{Proposition}
\newtheorem{theorem}[itlemma]{Theorem}
\newtheorem{itcorollary}[itlemma]{Corollary}
\newtheorem{itremark}[itlemma]{Remark}
\newtheorem{itremarks}[itlemma]{Remarks}
\newtheorem{itdefinition}[itlemma]{Definition}
\newtheorem{itexample}[itlemma]{Example}

\newenvironment{fact}{\begin{itfact}\rm}{\end{itfact}}
\newenvironment{claim}{\begin{itclaim}\rm}{\end{itclaim}}
\newenvironment{lemma}{\begin{itlemma}}{\end{itlemma}}
\newenvironment{remark}{\begin{itremark}\rm}{\end{itremark}}
\newenvironment{remarks}{\begin{itremarks} \rm}{\end{itremarks}}
\newenvironment{corollary}{\begin{itcorollary}}{\end{itcorollary}}
\newenvironment{proposition}{\begin{itproposition}}{\end{itproposition}}
\newenvironment{definition}{\begin{itdefinition}\rm}{\end{itdefinition}}
\newenvironment{example}{\begin{itexample}\rm}{\end{itexample}}
\newenvironment{proof}{\noindent {\em Proof}.\ \
}{\hspace*{\fill}$\halmos$\medskip}
\newcommand{\be}[1]{\addtocounter{for}{1} \begin{equation}\label{#1}}
\newcommand{\ee}{\end{equation}}
\newcommand{\bl}[1]{\begin{lemma}\label{#1}}
\newcommand{\br}[1]{\begin{remark}\label{#1}}
\newcommand{\brs}[1]{\begin{remarks}\label{#1}}
\newcommand{\bt}[1]{\begin{theorem}\label{#1}}
\newcommand{\bd}[1]{\begin{definition}\label{#1}}
\newcommand{\bp}[1]{\begin{proposition}\label{#1}}
\newcommand{\bc}[1]{\begin{corollary}\label{#1}}
\newcommand{\bfact}[1]{\begin{fact}\label{#1}}
\newcommand{\bex}[1]{\begin{example}\label{#1}}
\newcommand{\ec}{\end{corollary}}
\newcommand{\efact}{\end{fact}}
\newcommand{\eex}{\end{example}}
\newcommand{\el}{\end{lemma}}
\newcommand{\er}{\end{remark}}
\newcommand{\ers}{\end{remarks}}
\newcommand{\et}{\end{theorem}}
\newcommand{\ed}{\end{definition}}
\newcommand{\ep}{\end{proposition}}
\newcommand{\epr}{\end{proof}}
\newcommand{\bpr}{\begin{proof}}
\newcommand{\bcl}[1]{\begin{claim}\label{#1}}
\newcommand{\ecl}{\end{claim}}

\newcommand{\ecs}{\end{corollary}}
\newcommand{\eers}{\end{exercise}}
\newcommand{\eexs}{\end{example}}
\newcommand{\eems}{\end{example}}
\newcommand{\els}{\end{lemma}}
\newcommand{\eles}{\end{lemmaex}}
\newcommand{\ets}{\end{theorem}}
\newcommand{\eds}{\end{definition}}
\newcommand{\eps}{\end{proposition}}

\newcommand{\bi}{\begin{itemize}}
\newcommand{\ei}{\end{itemize}}
\newcommand{\ben}{\begin{enumerate}}
\newcommand{\een}{\end{enumerate}}

\def\vbar{\mathchoice{\vrule height6.3ptdepth-.5ptwidth.8pt\kern-.8pt}
   {\vrule height6.3ptdepth-.5ptwidth.8pt\kern-.8pt}
   {\vrule height4.1ptdepth-.35ptwidth.6pt\kern-.6pt}
   {\vrule height3.1ptdepth-.25ptwidth.5pt\kern-.5pt}}
\def\fudge{\mathchoice{}{}{\mkern.5mu}{\mkern.8mu}}
\def\bbc#1#2{{\rm \mkern#2mu\vbar\mkern-#2mu#1}}
\def\bbb#1{{\rm I\mkern-3.5mu #1}}
\def\bba#1#2{{\rm #1\mkern-#2mu\fudge #1}}
\def\bb#1{{\count4=`#1 \advance\count4by-64 \ifcase\count4\or\bba A{11.5}\or
   \bbb B\or\bbc C{5}\or\bbb D\or\bbb E\or\bbb F \or\bbc G{5}\or\bbb H\or
   \bbb I\or\bbc J{3}\or\bbb K\or\bbb L \or\bbb M\or\bbb N\or\bbc O{5} \or
   \bbb P\or\bbc Q{5}\or\bbb R\or\bbc S{4.2}\or\bba T{10.5}\or\bbc U{5}\or
   \bba V{12}\or\bba W{16.5}\or\bba X{11}\or\bba Y{11.7}\or\bba Z{7.5}\fi}}

\def \qed {{\hspace*{\fill}$\halmos$\medskip}}
\def \ra {\rightarrow }

\def \LL {{\cal{L}}}
\def \An {{\mathsf{A}}}
\def \A {{\cal{A}}}

\def \FF {{\cal{F}}}

\def \E {{\cal{E}}}
\def \C {{\cal{C}}}
\def \D {{\cal{D}}}
\def \M {{\cal{M}}}
\def \S {{\cal{S}}}

\newcommand{\ba}[1]{\addtocounter{for}{1} \begin{eqnarray}\label{#1}}
\newcommand{\ea}{\end{eqnarray}}

\def\sqr#1#2{{\vcenter{\vbox{\hrule height .#2pt
                             \hbox{\vrule width .#2pt height#1pt \kern#1pt
                                   \vrule width .#2pt}
                             \hrule height .#2pt}}}}

\def\pmb#1{\setbox0=\hbox{#1}%
   \kern-.025em\copy0\kern-\wd0
   \kern.05em\copy0\kern-\wd0
   \kern-.025em\raise.0433em\box0 }
\def\sqr#1#2{{\vcenter{\vbox{\hrule height.#2pt
     \hbox{\vrule width.#2pt height#1pt \kern#1pt
   \vrule width.#2pt}\hrule height.#2pt}}}}

\def\B{{\mathcal B}}
\def\N{{\mathbb N}}   
\def\Z{{\mathbb Z}}
\def\R{{\mathbb R}}
\def\Lb{{\mathbb L}_{b}}
\def\Li{{\mathbb L}}
\def\PP{{\mathbb P}}
  
\def\id{{\mathbf{1}}}

\def\e{\epsilon}                

\def\e{\epsilon}
\def\d{\delta}
\def\l{\lambda}
\def\lr{\lambda(\rho)}
\def\lnr{\lambda_n(\rho)}
\def\L{\Lambda}
\def\nur{\nu_{\rho}}

\def\g{\gamma}
\def\gr{\gamma_{\rho}}

\def\Gi{\Gamma_{\!\!\infty}}
\def\Gn{\Gamma_{\!\!n}}
\def\Gm{\Gamma_{\!\!m}}
\def\a{\alpha}

\def\v{\varphi}
\def\p{\partial}

\def\bs{\backslash}

\begin{document}

\title{On the Dirichlet problem for asymmetric zero range 
process on increasing domains.}
\author{Amine Asselah \\C.M.I., Universit\'e de Provence,\\
39 Rue Joliot-Curie, \\F-13453 Marseille cedex 13, France\\
asselah@cmi.univ-mrs.fr}
\date{}
\maketitle
\begin{abstract}
We characterize the principal eigenvalue of the 
generator of the asymmetric zero-range process on $\Z^d$ 
in dimensions $d\ge 3$, with Dirichlet boundary on special domains. 
We obtain a Donsker-Varadhan variational representation
for the principal eigenvalue, and
show that the corresponding eigenfunction is
unique in a natural class of functions. This allows
us to obtain asymptotic hitting time estimates.
\end{abstract}



{\em Running head}: Principal eigenvalue for asymmetric zero range.

\section{Introduction} 
We are concerned in this work
with obtaining hitting time estimates for the asymmetric zero-range 
process (AZRP). For this purpose, we consider the problem of
characterizing the principal eigenvalue and principal eigenfunctions
of the generator of AZRP, denoted by $\LL$, with Dirichlet boundary
on special domains. Though $\LL$ is 
neither compact, irreducible, nor self-adjoint, 
its physical origin endows crucial monotonicity properties. 

The AZRP models the conservative evolution of charged particles 
interacting over short range, in an electrical field. 
Thus, this process denoted by $\{\eta_t, t\ge 0\}$,
lives on $\{\eta:\eta(i)\in\N,
i\in \Z^d\}$, and evolves informally as follows.
At time zero and at each
site $i\in \Z^d$, we draw a number of particles $\eta(i)\in\N$.
To each particle we attach the trajectory of 
an asymmetric random walk with transition
kernel $\{p(i,j);i,j\in \Z^d\}$. Now, each site $i\in \Z^d$ has
an independent exponential process, its {\it clock},
of intensity $g(\eta_t(i))$ at time $t$, where $g:\N\to[0,\infty)$
is increasing.
When the clock of site $i$ rings, say at time $t$,
we choose a particle uniformely among the $\eta_t(i)$ ones and we move it
to its next position along its attached trajectory. 
The conservation of the particles number
imposes a one-parameter family of ergodic time-invariant measures $\{
\nur, \ \rho>0\}$, which happens to consist of product measures 
\cite{andjel,sethu}. The name zero-range is justified since only particles 
at the same site can interact with each other. Note also that $g(k)=k$
corresponds to independent random walks with clocks' intensity 1. 

A question motivated by physics is the time of occurrence of spots
with large densities of particles, say $\tau$,
when the gas is initially prepared with a homogeneous density. 
Thus, we consider a stationary process with respect to $\nur$,
and focus on occurrence time of patterns of the type
\be{eq0.5}
\A:=\{\eta: \sum_{i\in \S} \eta(i)>L\},\quad(\text{and}\quad
\tau:=\inf\{t: \eta_t\in \A\})
\ee
where the support of $\A$, $\S$, is a finite subset of $\Z^d$, and $L$ a given integer.

The key feature of our model is that the partial order --$\eta\prec \zeta$
meaning $\eta(i)\le \zeta(i)$ for all $i\in \Z^d$-- is preserved
under the evolution. Another important feature is
that the invariant measures $\{\nur,\ \rho>0\}$ 
all satisfy FKG's inequality, i.e.
for $f$ and $g$ increasing functions
\be{def-FKG}
\int fg d\nur\ge \int fd\nur \int g d\nur.
\ee
This was the setting of \cite{ac} whose relevant results we now recall. 
A simple subadditive argument yielded the asymptotic rough estimate
\be{eq0.6}
\lambda(\rho)=-\lim_{t\to\infty} \frac{1}{t}\log(P_{\nur}(\tau>t)).
\ee
When the drift is nonzero, $\lambda(\rho)$ is positive in
any dimensions. Furthermore,
if we denote by $\LL^*$ the dual of $\LL$ in $L^2(\nur)$, which corresponds
to an AZRP with reversed drift, then when dimension $d\ge 3$,
there exist $u,u^*\in L^p(\nur)$ for any $p\ge 1$ in the domain of
$\LL$ and $\LL^*$ respectively, with
\be{eq0.7}
\text{(i)}\quad 1_{\A^c}\LL(u)+\lambda(\rho)u=0,\qquad\text{and}\qquad
\text{(ii)}\quad 1_{\A^c}\LL^*(u^*)+\lambda(\rho) u^*=0,
\ee

However, and this was most unfortunate from a physical point of view,
a link with finite dimensional dynamics was missing, as well as
a variational representation for $\lambda(\rho)$. 
This is what we establish in this paper. 
Moreover, we establish uniqueness for $u$ in some class of functions,
which in turn yields an asymptotic estimate for the hitting time.

We have chosen to introduce
some symbols intuitively so as to be able to state our main results 
postponing definitions and notations as much as possible
to Section~\ref{sec-not}.

A way of defining the AZRP with initial law $\nur$
on $\Z^d$ is through a limit of
irreducible processes, where particles evolve on $[-n,n]^d$ as
a zero-range process with creation and annihilation at the boundary.
Informally, if $\FF_n$ is the $\sigma$-field generated by
$\{\eta(i), i\in [-n,n]^d\}$, then we define
\[
\LL_n^{\rho}(\v)=E_{\nur}[\LL(\v)|\FF_n].
\]
The generator $\LL_n^{\rho}$ will be shown to inherit the same property
of monotonicity as $\LL$ and to have $\nur$ as invariant measure.
Thus, its principal Dirichlet eigenvalue $\lnr$ is obtained
as in (\ref{eq0.6}). We show in Section~\ref{varia-ln}
that $\LL_n^{\rho}$ has a unique normalized eigenfunction $u_n\ge 0$,
associated with $\lnr$.
Then, our main observation in Section~\ref{approx-ln} is the following.
\bl{lem-approx} For $\lambda(\rho)$ given by (\ref{eq0.6}),
and $\lnr$ corresponding to $\LL_n^{\rho}$, we have
\be{eq2.3} 
\lim_{n\to\infty} \lnr=\inf_n \lambda_n(\rho)=\lambda(\rho).
\ee
\el
Moreover, we establish a link between finite and infinite volume
eigenfunctions.
\bt{the6} When $d\ge 3$, $\{u_n,n\in\N\}$ converges to a 
solution of (\ref{eq0.7}(i)) in weak-$L^2(\nur)$.
\et
In \cite{ac}, a solution of (\ref{eq0.7}(i)) was obtained through
another sequence, say $\{u_t,t\ge 0\}$ in which $u_t$ was the density
(w.r.t $\nur$) of the law of time-reversed process $\eta_t^*$ conditioned on
$\{\tau>t\}$. The functions
$\{u_t,t\ge 0\}$ where positive and decreasing on $\A^c$, and satisfied
the following uniform bound: for site $i$ large enough,
if $\e_i$ is the probability that a random walk starting on $i$ with 
transition kernel $\{p(.,.)\}$ hits the support of $\A$, then when $d\ge 3$
\be{def-D}
0\le u_t(\eta)-u_t(A_i^+\eta)\le \e_i u_t(\eta),
\ee
where $A_i^+$ adds a particle at $i\in\Z^d$.

We denote by $\D_{\rho}$ the convex
set of non-negative decreasing functions of finite integral (w.r.t $\nur$),
satisfying (\ref{def-D}). We denote by $\D^+_{\rho}$ the positive
functions of $\D_{\rho}$. Finally, we define
a {\it dual} space of probability measures, $\M_{\rho}$,
absolutely continuous with respect to $\nur$, and whose density
satisfies a condition similar to (\ref{def-D}). 

Intuitively, a Donsker-Varadhan's type functional
would read $\Gi(\v,\mu)=$``$\int \LL(\v)/\v d\mu$'' for $(\v,\mu)\in
\D^+_{\rho}\times \M_{\rho}$. One problem is that $\LL$ cannot
be defined on $\D^+_{\rho}$ as a convergent series. Thus,
we define $\Gi(\v,\mu)$ in
Proposition~\ref{prop4.1} of Section~\ref{sec-dv}
as a Cauchy limit using cancelation due to gradient bounds (\ref{def-D})
on $\v$ and $d\mu/d\nur$. 

We obtain in Section~\ref{sec-the2} a 
Donsker-Varadhan variational formula for the principal eigenvalue.
\bt{the2} When $d\ge 3$, and $\A$ is increasing with bounded support, we have
\be{eq0.15}
\lambda(\rho)=-\sup_{\mu\in \M_{\rho}}\inf_{\v\in \D^{+}_{\rho}}
\Gi(\v,\mu).
\ee
\et
Obtaining (\ref{eq0.15}) is linked with the issue of uniqueness of
the principal eigenfunction, since the minimax theorem hidden
behind Donsker-Varadhan formula requires a convex functional
$h\mapsto \Gi(e^{h},\mu)$, on a convex set of functions regular enough.
Note that $\D^+_{\rho}$ is all the more appropriate since when written
for $h=\log(\v)$ with $\v\in\D^+_{\rho}$, condition (\ref{def-D}) reads
\be{def-H}
h(\eta)\ge h(\An_i^+\eta)\ge h(\eta)+\log(1-\e_i)\quad(\text{when }
\e_i<1),
\ee
and define a convex set. Now, the main uniqueness result is the following.
\bt{the3} When $d\ge 3$, there is a unique normalized
Dirichlet eigenfunction in $\D_{\rho}$. This eigenfunction
is positive $\nur$-a.s.\ on $\A^c$.
\et
The proofs of Theorem~\ref{the3} and Theorem~\ref{the6}
are conducted in Section~\ref{unique}.
We sketch the simple intuitive steps behind the proof of uniqueness.
Assume there exist $u,\tilde u$ solutions of (\ref{eq0.7}(i)) in
$\D_{\rho}$. Then, they are actually positive (on $\A^c$), and satisfy
\be{sketch2}
\forall\mu\in\M_{\rho},
\qquad\Gi(u,\mu)=\Gi(\tilde u,\mu)=-\lr.
\ee
As already mentionned, if $u,\tilde u\in \D_{\rho}$
and $\gamma\in]0,1[$, then
$u_{\gamma}:=u^{\gamma}\tilde{u}^{1-\gamma}\in \D_{\rho}$. Now,
by convexity of $h\mapsto \Gi(\exp(h),\mu)$
\be{sketch3}
\forall \mu\in\M_{\rho},\qquad
-\lr=\gamma \Gi(u,\mu)+(1-\gamma)\Gi(\tilde u,\mu)\ge
\Gi(u_{\gamma},\mu)
\ee
We now choose a special $\mu$ so that equality obtains in (\ref{sketch3}).
The space $\M_{\rho}$ is built so that if $u^*$ is a positive solution
of (\ref{eq0.7}(ii)), then
\be{sketch4}
d\mu^*:=\frac{u_{\gamma}u^*}{\int\!\!u_{\gamma}u^*d\nur}d\nur\in
\M_{\rho}.
\ee
Now, by a formal use of duality
\be{sketch5}
\Gi(u_{\gamma},\mu^*)=
``\int \frac{\LL(u_{\gamma})}{u_{\gamma}}
\frac{u_{\gamma}u^*}{\int\!\!u_{\gamma}u^*d\nur}d\nur\text{''}=
``\int \frac{\LL^*(u^*)}{u^*}
\frac{u_{\gamma}u^*}{\int\!\!u_{\gamma}u^*d\nur}d\nur\text{''} =-\lr.
\ee
Finally, the case of equality in (\ref{sketch3}) implies that
$\tilde u/u$ is $\nur$-a.s.\ constant on $\A^c$ by using the triviality
of the $\sigma$-field of exchangeable events under $\nur$.

As a consequence of Theorem~\ref{the3}, we obtain an asymptotic estimate
of the first hitting time of $\A$. To link this last result with
those of~\cite{ac}, we recall Corollary 2.8 of~\cite{ac} which
was based on $L^p(\nur)$ estimates for $u$ and $u^*$. When
$d\ge 3$,
there is a positive constant $\underline{c}$ such that for any $t\ge 0$,
\be{eq0.16}
\underline{c}\le \exp(\lr t) P_{\nur}(\tau>t)\le 1.
\ee
As a corollary of the uniqueness of the principal eigenfunction in $\D_{\rho}$,
we obtain the following estimates whose proof makes up Section~\ref{hitting}.
\bt{the4} When $d\ge 3$, 
\be{eq0.17}
\lim_{t\to\infty} \frac{1}{t}\int_0^t e^{\lr s} P_{\nur}(\tau>s)ds
=\frac{1}{\int\!\!u u^* d\nur}.
\ee
\et
\section{Notations and preliminaries.}
\label{sec-not}
We first recall in Section~\ref{sec-zero}, the hypotheses needed
to define the AZRP on $\Z^d$. Then, in Section~\ref{sec-pattern},
we describe the class of patterns we consider here.
Section~\ref{sec-function} contains the definition of all function
spaces which we use.
\subsection{The zero-range process}
\label{sec-zero}
The transition kernel $\{p(i,j),\ i,j\in \Z^d\}$ is associated with 
a single-particle trajectory and satisfies for all $i,j$ in $\Z^d$
\ba{def-p}
\text{(i)}&& p(i,j)\geq 0,\quad p(i,i)=0,\quad 
\textstyle{\sum_{i\in \Z^d}} p(0,i)=1. \cr
\text{(ii)}&&p(i,j)=p(0,j-i)\quad\text{(translation invariance)}.\cr
\text{(iii)}&& p(i,j)=0\text{ if }|i-j|>R\text{ for some fixed }R
\quad\text{(finite range)}.\cr
\text{(iv)}&& {\rm If\ }p_s(i,j)=p(i,j)+p(j,i),\text{ then}\quad
\forall i\in \Z^d,\ \exists n,\quad p_s^{(n)}(0,i)>0
\quad\text{(irreducibility)}.\cr
\text{(v)}&& \textstyle{\sum_{i\in \Z^d}}
ip(0,i)\not= 0 \quad\text{(positive drift)}.
\ea
Note that by (i) and (ii), the transition kernel $p(.,.)$ is
doubly stochastic. Thus, we can introduce a {\it dual} transition kernel
$\{p^*(i,j),\ i,j\in \Z^d\}$, with $p^*(i,j)=p(j,i)$.

We also need a particle dependent intensity $g$ which satisfies 
\ba{def-g}
\text{(i)}&&g:\N\to [0,\infty) \text{ is increasing}.\cr
\text{(ii)}&&g(0)=0, \quad g(1)=1\quad\text{(normalization)}.\cr
\text{(iii)}&& \quad \Delta:=\sup_k\left(g(k+1)-g(k)\right)<\infty.
\ea
For notational simplicity, we call the intensity
at site $i\in \Z^d$, $g_i(\eta):=g(\eta(i))$.

For any $\g\in [0,\sup_k g(k)[$, we define
a probability $\theta_{\g}$  on $\N$, by
\be{def-m}
\theta_{\g}(0)=1/Z(\g),\quad\text{ and when }n\not= 0,\quad
\theta_{\g}(n)=\frac{1}{Z(\g)}\frac{\g^n}{g(1)\dots g(n)},
\ee
where $Z(\g)$ is the normalizing factor. If we set
$\rho(\g):=\sum_{n=1}^{\infty} n\theta_{\g}(n)$, then
$\rho:[0,\sup_k g(k)[\to [0,\infty[$ is increasing.
Let $\gamma(.)$
be the inverse of $\rho(.)$, and for a constant density $\rho>0$,
let $\nur$ be the product
probability with marginal law $\theta_{\g(\rho)}$. Thus, we have
\be{def-nur}
\forall B\subset \Z^d,\quad\int\prod_{i\in B}
\eta(i)d\nur=\rho^{|B|},\quad{\rm and}\quad
\int\!\! g_i(\eta) \v(\An_i^-\eta)d\nur(\eta)=\g(\rho)
\int\!\! \v d\nur,
\ee
where $\An_i^-\eta$ has one particle less than $\eta$ at site $i$.
Also, we will often use that
\be{eq-obvious}
0\le g(n)\le \Delta n,\quad(\text{by (ii) and (iii) of (\ref{def-g})}),\qquad
\text{and}\qquad \int\!\! g_i^p d\nur<\infty,\quad\text{for any }p\in\N.
\ee
Following \cite{liggett0},
(see also \cite{andjel} and \cite{sethu} Section 2), let
\[
\alpha(i)=\sum_{n=0}^{\infty} 2^{-n} p^n(i,0),\quad
\text{and for }\eta,\zeta\in \N^{\Z^d},
\quad ||\eta-\zeta||=\sum_{i\in \Z^d}|\eta(i)-\zeta(i)|
\alpha(i).
\]
Since the transition kernel $p$ is finite range (by \ref{def-p}(iii)), 
another possible
choice is $\alpha(k)=\exp(-(|k_1|+\dots+|k_d|))$ for any 
site $k=(k_1,\dots,k_d)$ (see \cite{liggett0}).
Our state space is $\Omega=\{\eta:||\eta||<\infty\}$, and we call
$\Li$ the space of Lipshitz functions from $(\Omega,||.||)$
to $(\R,|.|)$, and $\Lb$ the subspace of $\Li$ consisting
of bounded functions. For $\v\in\Li$, we call
\be{def-Lip}
L(\v):=\sup\{\frac{|\v(\eta)-\v(\xi)|}{||\eta-\xi||}:
||\eta-\xi||>0,\ \eta,\xi\in \Omega\}.
\ee
In \cite{andjel}, it is shown that a semi-group
can be constructed on $\Li$ with formal generator
\be{def-formal}
\LL \v(\eta):=\sum_{i,j\in \Z^d}p(i,j)g(\eta(i))\left(
\v(T^i_j\eta)-\v(\eta)\right),
\ee
where $T^i_j\eta(k)=\eta(k)$ if $k\not\in \{i,j\}$, 
$T^i_j\eta(i)=\eta(i)-1$, and $T^i_j\eta(j)=\eta(j)+1$.
If we set $\nabla^i_j \v=\v\circ T^i_j-\v$, we
will often use that on $\{\eta(i)>0\}$
\be{link-nabla}
\nabla^i_j \v=(\v\circ\An_j^+-\v\circ \An_i^+)
\circ\An_i^-.
\ee
Thus, if we set $\Delta^j_i\v=\v\circ\An_j^+-\v\circ \An_i^+$, 
and use (\ref{eq-obvious}) and (\ref{link-nabla}),
we have the following integration by parts formula
\be{eq-byparts}
\int\!\! g_i \nabla^i_j(\v)f d\nur=\gr \int\!\! \Delta^j_i(\v)\ 
\An_i^+(f) d\nur.
\ee
Also, for convenience, we often write $\An_i^{\pm}\v$ for $\v\circ\An_i^{\pm}$.

In \cite{sethu} Section 2, $\LL$ is extended to a generator, 
again called $\LL$ for convenience,
on $L^2(\nur)$ for any $\rho>0$. It is also shown that $\Lb$ is
a core for $\LL$. Moreover, $\{\nur, \rho>0\}$ 
are ergodic invariant measures for $\LL$. We denote by
$\D(\LL,L^2(\nur))$ the domain of $\LL$ in $L^2(\nur)$, and
by $||.||_{\nu}$ the $L^2(\nu)$-norm, for any probability measure
$\nu$. Finally, we consider the adjoint (or time-reversed)
of $\LL$ in $L^2(\nur)$, acting on Lipshitz
functions $\v$ and $\psi$ by
\be{def-adjoint}
\int \LL^*(\v) \psi d\nur:=\int \v\LL(\psi)d\nur.
\ee
With our hypothesis, $\LL^*$ is again the generator of a 
zero-range process with transition kernel $p^*(.,.)$ 
satisfying $p^*(i,j):=p(j,i)$ and with the same function $g$. We denote
by $\{S_t^*\}$ the associated semi-group, and by $P^*_{\eta}$ the 
associated Markov process with initial configuration $\eta\in \Omega$.

\subsection{Special patterns.}
\label{sec-pattern}
We first recall that there is a partial order on $\Omega$.
For $\eta,\xi \in \Omega$, we say that $\eta \prec \xi$
if $\eta(i) \leq \xi(i)$ for all $i \in \Z^d$. A
function $f:\Omega\to \R$ is increasing if for $\eta \prec \xi$,
$f(\eta)\leq f(\xi)$. Also, we say that $A \subset \Omega$ is
increasing if its indicator $1_A$ is increasing. Finally, for
given probability measures $\nu,\mu$ on $\Omega$, we say that
$\nu \prec \mu$ if $\int f d\nu \leq \int f d\mu$ for every
increasing function $f$. The zero-range process
is a monotone process, i.e. 
there is a coupling such that $P_{\eta,\zeta}(\eta_t\le \zeta_t,
\forall t)=1$ whenever $\eta\le \zeta$.

We will be concerned with the hitting time of pattern, $\A$, with
the following properties dubbed $(\C\!\!-\!\!\FF)$ for connectedness and
finiteness:
\begin{itemize}
\item[(i)] It is non-empty, and its support $\S$ is bounded. Thus,
$\nur(\A)>0$.
\item[(ii)] It is increasing, and $0_{\S}:=\{\eta: \eta(i)=0,
\forall i\in \S\}\not\subset \A$. Thus, $\nur(\A)<1$.
\item[(iii)] Its complement, $\A^c$, is connected, and is partitioned
into a finite number of cylinders with support in $\S$, whose
set we denote by $\Theta$. In other words, for any
cylinder $\theta\in\Theta$, there is an integer $n$, a sequence
$\theta_0,\dots,\theta_n\in \Theta$, and $i_1,\dots,i_n\in \S$ such
that
\[
\theta_0:=0_{\S},\quad \theta_n=\theta,\quad\text{and}\quad
\theta_k=\An_{i_k}^+\theta_{k-1}\quad\text{for }k=1,\dots,n.
\]
\end{itemize}
A typical example of patterns satisfying $(\C\!\!-\!\!\FF)$ is given
in (\ref{eq0.5}). Note also that if $\A$ satisfies $(\C\!\!-\!\!\FF)$,
there is an integer $L$ such that $\{\eta:\sum_{\S}\eta(i)>L\}\subset
\A$.

We denote by $\bar \LL:=1_{\A^c}\LL$ and
$\{\bar S_t,t\ge 0\}$, respectively the generator and associated semi-group
for the process killed on $\A$. 
\subsection{Function spaces.}
\label{sec-function}
The topology on $\{\eta:\eta(i)\in\N, i\in \Z^d\}$, is the product
of discrete topology, so that $\{\eta_n,n\in\N\}$ converges
to $\eta$, if for any site $i\in \Z^d$, there is $n_0$ such that
for $n\ge n_0$ $\eta_n(i)=\eta(i)$.

Let $H_{\S}:=\inf \{t: X_t\in \S\}$ for $\{X_t\}$ a random
walk with transition kernel $\{p(i,j);i,j\in \Z^d\}$. Note that
$\e_i:=\PP_i(H_{\S}< \infty)\to 0$ as $||i||\to\infty$, 
(as well as $\e^*_i$ corresponding to a reversed drift) and when
the dimension $d\ge 3$, then we have the classical results
\[
\sum_{i\in \Z^d} \e_i^2+(\e_i^*)^2<\infty.
\]
Let $\A$ satisfy $(\C\!\!-\!\!\FF)$.
Choose $n$ large enough so that $\S\subset \L_n:=[-n,n]^d$, and set 
$\Omega_n=\{\eta:\L_n\to \N\}$, and $\FF_n:=\sigma(\{\eta(i), i\in \L_n\})$. 
We often make the abuse of considering functions on $\Omega_n$ as defined
also on $\Omega_m$ for $m\ge n$, but depending only on the
sites of $\L_n$.
\subsubsection{Functions on $\Omega_n$.}
A function $\v$ on $\L_n$ with $\v|_{\A}\equiv 0$ belongs to $\D_n$ when
\ba{def-Dn}
\text{(0)}&& 0\le \v,\cr
\text{(i)}&&\forall \eta,\zeta\in\Omega_n\bs\A,\text{ if }\eta\prec \zeta
\quad\text{ then}\quad \v(\zeta)\le \v(\eta),\cr
\text{(ii)}&&\forall \eta\in\Omega_n\bs\A,\ \forall i\in\L_n\bs\S,
\qquad\v(\eta)-\v(\An_i^+\eta)\le \v(\eta) \e_i,\cr
\text{(iii)}&&\int \v d\nur<\infty.
\ea
When $\e_i^*$ replaces $\e_i$ in (ii), we say that $\v$ belong to $\D_n^*$.
Also, we set $\D_n^+:=\D_n\cap\{\v$ positive on $\A^c\}$. 
\bl{lem3.1} $\D_n$ is a convex subset of $\Lb$.
When $d\ge 3$, if $\v\in \D_n^+$, then $\v$ and 
$1_{\A^c}/\v$ are in $L^p(\nur)$ for any $p\ge 1$.
\el
\bpr
If $\v\in \D_n$, note that $\v$ is bounded since $0\le \v(\eta)\le
\v(0_{\L_n})$, where $0_{\L_n}$ is the empty configuration of $\Omega_n$.
Take $\eta,\zeta\in \Omega_n\bs\A$, 
and let $\xi=\eta\vee \zeta-\eta\wedge \zeta$,
and set $m=\sum_i \xi(i)$. Since $\v$ is decreasing
\[
|\v(\eta)-\v(\zeta)|\le \v(\eta\wedge \zeta)-\v(\eta\vee \zeta).
\]
Now, let $\{\eta_i,i=0,\dots,m\}$ be an ordered sequence with
\[
\eta\wedge \zeta=\eta_0\prec \eta_1\prec\dots\prec \eta_m=\eta\vee \zeta,
\quad\text{with}\quad \eta_i=\An_{j_i}^+\eta_{i-1},
\]
where $\{j_i,i=1,\dots,m\}$ are the positions of the $m$ particles of $\xi$. Then,
\[
\v(\eta\wedge \zeta)-\v(\eta\vee \zeta)\le \sum_{i=0}^{n-1}
\v(\eta_i)-\v(\eta_{i+1})\le \sum_{i=1}^{n}\v(\eta_{i-1})\e_{j_i}.
\]
We use that $\v(\eta_i)\le \v(0_{\L_n})$, and that
$\sum_i \e_{j_i}=\sum_k \e_k \xi(k)$. Thus,
\be{eq-lip}
|\v(\eta)-\v(\zeta)|\le \v(0_{\L_n})\sum_{k\in \L_n} \e_k \xi(k)\le
\v(0_{\L_n})\sup_{k\in \L_n}(\frac{\e_k}{\alpha_k})\ 
\sum_{k\in \L_n}\xi(k) \alpha(k).
\ee
Now, if $\eta,\zeta\in \A$, then (\ref{eq-lip}) holds.
Assume that $\eta\in \Omega_n\bs\A$ but $\zeta\in\A$. Inequality
(\ref{eq-lip}) follows once we notice that $||\eta-\zeta||\ge
\inf_{\S} \a>0$. Thus, $\v$ is a Lipshitz bounded function.
Now, $\v$ and $1_{\A^c}/\v$ are in $L^p(\nur)$ for any
integer $p$ by Lemmas~\ref{lem5.1} and~\ref{app-key} of the Appendix.
\epr

For any $\v\in \D^{+}_n$, we can define its logarithm on $\A^c$, 
$h=\log(\v)$; on $\A$ we set $h\equiv-\infty$. 
Note that (\ref{def-Dn}) reads for $h$
\ba{def-En}
\text{(i)}&&\forall \eta,\zeta\in\Omega_n\bs\A,\text{ if }\eta\prec \zeta
\quad\text{ then}\quad h(\zeta)\le h(\eta),\cr
\text{(ii)}&&\forall \eta\in\Omega_n\bs\A,\ \forall i\in\L_n\bs\S,
\qquad h(\An_i^+\eta)\ge h(\eta) +\log(1-\e_i),\cr
\text{(iii)}&&\int \exp(h) d\nur<\infty.
\ea

Thus, we will say that $h\in\E_n$ if it satisfies (\ref{def-En}).
A key and simple observation is the following.
\bl{lem3.2} $\E_n$ is a convex set.
\el
\bpr
Inequalities (\ref{def-En}) (i) and (ii)
are stable under convex combination. Also,
for $\gamma\in ]0,1[$, and $h_1,h_2\in \E_n$ by H\"older inequality
\be{eq3.3}
\int \exp(\gamma h_1+(1-\gamma) h_2)d\nur\le
\left(\int e^{h_1}d\nur\right)^{\gamma}
\left(\int e^{h_2}d\nur\right)^{1-\gamma}<\infty.
\ee
\epr

We now define $\M_n$ a space of probability measures
whose elements have a density with respect to $\nur$,
generically noted $f$ satisfying: (i) $f$ is decreasing on $\A^c$,
$f|_{\A}\equiv 0$, and 
\be{def-Mn}
\text{(ii)}\qquad
\forall \eta\in\Omega_n\bs\A,\ \forall i\not\in\S
\qquad f(\eta)-f(\An_i^+\eta)\le f(\eta) (\e_i+\e_i^*)
\ee
\bl{lem3.6} Assume that $d\ge 3$.
$\M_n$ is a convex and compact set in the weak topology.
\el
\bpr
The convexity of $\M_n$ is obvious. 
Consider the compact decreasing set
\be{def-K}
K_M=\{\eta\in \Omega_n: \eta(i)\le M,\ \forall i\in \L_n\}.
\ee
Note that $\M_n$ is tight:
\[
\lim_{M\to\infty} \sup_{\mu\in \M_n} \mu(K_M^c) =0.
\]
Indeed, since $d\mu/d\nur$ is decreasing for any $\mu\in\M_n$,
by FKG's inequality
\[
\forall \mu\in\M_n,\qquad
\mu(K_M^c)=\int \id_{K_M^c}\frac{d\mu}{d\nur}d\nur\le 
\nur(K_M^c)
\stackrel{\scriptstyle{M \ra \infty}}{\longrightarrow}0.
\]
Let $\{\mu_n,n\in\N\}$ be in $\M_n$, with densities
$\{f_n:=d\mu_n/d\nur\}$. Let $\{\mu_{n_k}\}$ a converging
subsequence to $\mu$. For any $\eta\in \Omega_n$, $\id_{\eta}$ is a 
bounded continuous function, so that
\be{eq3.22}
f_{n_k}(\eta)\nur(\eta)=\int \id_{\eta}d\mu_{n_k}
\stackrel{\scriptstyle{k \ra \infty}}{\longrightarrow}
\mu(\eta)=f(\eta)\nur(\eta).
\ee
Thus, $f_{n_k}$ converges pointwise to $f$ on $\Omega_n$. It is clear that
$f$ satisfies (\ref{def-Mn}) so that $\mu\in\M_n$.
\epr

An important feature of $\M_n$ is the following.
\bl{lem3.3} Assume that $d\ge 3$. 
If $\v\in \D^{+}_n$ and $\v^*\in (\D^{*}_n)^+$, then
\be{eq3.4}
d\mu=\frac{\v\v^*d\nur}{\int \v\v^*d\nur}\in \M_n.
\ee
\el
\bpr
First, by Lemma~\ref{lem5.1}, $\int\!\!\v \v^*d\nur<\infty$. Also, 
note that $\v,\v^*>0$ on $\A^c$ so that $\int \v\v^*d\nur>0$. Thus,
$\mu$ given in (\ref{eq3.4}) is well defined.
Now, since $\v$ and $\v^*$ are decreasing on $\A^c$ and
positive, $d\mu/d\nur$ is decreasing on $\A^c$. Now, if $\zeta=\An_i^+\eta$,
for $i\not\in \S$
\ba{eq3.5}
\v(\eta)\v^*(\eta)-\v(\zeta)\v^*(\zeta)&=&
\v^*(\eta)(\v(\eta)-\v(\zeta))+ \v(\zeta)(\v^*(\eta)-\v^*(\zeta))\cr
&\le&  \v(\eta)\v^*(\eta)(\e_i+\e_i^*).
\ea
Thus, $\mu$ satisfies (i) and (ii) of (\ref{def-Mn}). 
\epr

\subsubsection{Functions on $\Omega$.}
We define $\D_{\rho}$ as the natural extention of $\D_n$ to functions
defined on the whole of $\Omega$. Thus, functions in 
$\D_{\rho}$ satisfy the inequalities in
(\ref{def-Dn}(0)-(iii)) but almost surely with respect to $\nur$. Also,
$\D^+_{\rho}$ denotes the functions of $\D_{\rho}$ positive
$\nur$-a.s.\ on $\A^c$. Similarly, we extend $\M_n$ into $\M_{\rho}$,
the space of probability measures absolutely continuous with respect
to $\nur$, whose densities satisfy $\nur$-a.s.\ the same
conditions as function of $\M_n$, but extended on the whole of $\Z^d$.
Note that by linearity of the conditional expectation,
for $\v\in \D_{\rho}$, $E_{\nur}[\v|\FF_n]\in \D_n$, and
similarly if $\mu\in \M_{\rho}$ with density $f$, then
$E_{\nur}[f|\FF_n]d\nur\in \M_n$.
\bl{lem-compact}
$\M_{\rho}$ is compact in the weak topology.
\el
\bpr
First, by Remark~\ref{rem-M} of the Appendix,
there is a constant $C(\rho,2)>0$ such that
\[
\sup_{\mu\in \M_{\rho}} \int\!\!(\frac{d\mu}{d\nur})^2 d\nur\le
C(\rho,2).
\]
Recall that by Banach-Alaoglu Theorem, $\{d\mu/d\nur,\ \mu\in \M_{\rho}\}$
is weak-$L^2(\nur)$ compact in $L^2(\nur)$. 
Secondly, for any $\mu\in \M_{\rho}$ and integer $n$, as already
mentionned
\[
d\mu^{(n)}:=E_{\nur}[\frac{d\mu}{d\nur}\big| \FF_{n}]d\nur\in \M_n.
\]

Now, let $\{\mu_k,k\in \N\}$ be in $\M_{\rho}$, and
let $\mu_{\infty}$ be a weak-$L^2(\nur)$ limit along a subsequence,
say $\{ n_k\}$. Note that for each integer $n$, the following
convergence holds in weak-$L^2(\nur)$
\[
f_{n_k}^{(n)}:=E_{\nur}[\frac{d\mu_{n_k}}{d\nur}\big| \FF_{n}]
\stackrel{\scriptstyle{k \ra \infty}}{\longrightarrow}
f_{\infty}^{(n)}:=E_{\nur}[\frac{d\mu_{\infty}}{d\nur}\big| \FF_{n}].
\]
Moreover, $f_{\infty}^{(n)}d\nur\in \M_n$, since $\M_n$
is compact by Lemma~\ref{lem3.6}. Finally, the sequence
$\{f_{\infty}^{(n)},\ n\in \N\}$ is a positive martingale which,
by the martingale convergence Theorem, converges $\nur$-a.s.\ to
$f_{\infty}$. Clearly, inequality (\ref{def-Mn}) holds $\nur$-a.s.\ for
$f_{\infty}$.
\epr
\br{rem-comp}
With the same arguments, we obtain that $\D_{\rho}\cap\{\v:\int\!\!\v
d\nur\le c\}$ is weak-$L^2(\nur)$ compact, for any constant $c>0$.
\er
\br{rem-existence}
We give now more details on how a solution $u$ to
(\ref{eq0.7}(i)) was obtained in~\cite{ac}, 
and why $u\in \D_{\rho}$ actually.
We recall that for any probability $\mu$,
$\Phi(\mu)$ introduced in~\cite{fkmp} was
the invariant measure of the renewal process corresponding
to $\{\eta_t\}$ started afresh from measure $\mu$
each time it hits $\A$. Also, for any integer $k$, the map $\Phi^{(k)}$ was
the $k$-th iterates of $\Phi$. It is shown in Theorem 2.4
of \cite{ac} that the Cesaro weak-$L^2(\nur)$
limits of $\{\Phi^{(k)}(\nur),k\in \N\}$ are solutions of (\ref{eq0.7}(i)).
There is actually a simple expression
for $\Phi^{(k)}$. Since $\lambda(\rho)>0$, we have
$\int_0^{\infty} P_{\nur}(\tau>t) t^k dt<\infty$, and
the following probability $dm_k(t)$ on $\{t\ge 0\}$ is well defined
\be{eq3.9}
dm_k(t)=\frac{P_{\nur}(\tau>t) t^k dt}{\int_0^{\infty}
P_{\nur}(\tau>t) t^k dt}\qquad\text{and}\qquad
\frac{d\Phi^{(k)}(\nur)}{d\nur}(\eta)=\int_0^{\infty}
u_{t}(\eta) dm_k(t),
\ee
where $u_t$ is mentionned in the paragraph preceding (\ref{def-D}).
Since, $u_{t}\in \D_{\rho}$, it is clear that for any integer $k$,
$d\Phi^{(k)}(\nur)/d\nur\in\D_{\rho}$ as well as the Cesaro mean
since $\D_{\rho}$ is convex. 
Now, since $\Phi^{(k)}(\nur)$ are probability measures,
Remark~\ref{rem-comp} implies that all
the Cesaro limits are in $\D_{\rho}$. Thus, there exists
a solution of (\ref{eq0.7}(i)) in $\D_{\rho}$: we denote it by $u$.
Notice also that our uniqueness result, Theorem~\ref{the3}, implies
that the whole Cesaro limit converges to $u$, thus strengthening
the results of \cite{ac}.
\er
\section{From finite domains to $\Z^d$.}
\subsection{Irreducible dynamics on $\L_n$.}
\label{irreducible}
Following the approach of~\cite{lig-spi}, as in~\cite{andjel},
we first consider, for any integer $k$ and $m$,
a finite-state generator $\LL^k_{(m)}$ on the
hyper-surface
\[
\Omega^k_{(m)}:=\{ \eta\in \N^{\L_m}:\ \sum_{i\in \L_m}
\eta(i)=k\}.
\]
For this purpose we introduce, for any integer $n$ and for $i,j\in
\L_n$
\be{def-pn}
p_n(i,j):=
\left\lbrace\begin{array}{l}
p(i,j)\quad\text{if}\quad i\not= j\\
\sum_{k\not\in\L_n} p(i,k)\quad\text{if}\quad i= j
\end{array}\right.,\quad\text{and}\quad
p_n^*(i,j):=
\left\lbrace\begin{array}{l}
p^*(i,j)\quad\text{if}\quad i\not= j\\
\sum_{k\not\in\L_n} p^*(i,k)\quad\text{if}\quad i= j.
\end{array}\right.
\ee
Note that $\{p_n(i,j)\}$ is not doubly stochastic.
The $\{\LL^k_{(m)}, k\in \N\}$ have the same expression,
though on different domains 
\be{def-Llig}
\forall \eta\in \Omega^k_{(m)},\qquad
\LL^k_{(m)}(\v)(\eta)= 
\sum_{i,j\in \L_m}p_m(i,j) g_i(\eta)(\v(T^i_j\eta) -\v(\eta)).
\ee
The process generated by $\LL^k_{(m)}$ is well defined.
Now, we take $n<m-R$, where $R$ is the range of the transition
kernel $p(.,.)$, and for $\v\in \D_n$, we define
\be{def-Lamin}
\LL_n^{\rho}(\v)=\lim_{K\to\infty}\sum_{k=0}^K
E_{\nur}[1_{\Omega^k_{(m)}}\LL^k_{(m)}(\v)|\FF_{\L_n}].
\ee
This limit is well define since $\D_n\subset \Lb$, and
\[
p_m(i,j) g_i(\eta)(\v(T^i_j\eta) -\v(\eta))\le
L(\v)p_m(i,j) g_i(\eta)(\a(i)+\a(j)),
\]
so that by Lemma 2.1 of \cite{sethu}, we have that
\[
\sum_{k\ge 0} \int \left(\LL^k_{(m)}(\v)\right)^2 1_{\Omega^k_{(m)}}
d\nur<\infty.
\]
Also, the expression $\LL^k_{(m)}(\v)$, and
the limit (\ref{def-Lamin}) are independent of $m$ when $m> n+R$,
and we called the latter $E_{\nur}[\LL(\v)|\FF_{\L_n}]$ in the Introduction.
Since $\{\LL^k_{(m)}, k\in \N\}$ have the same expression, 
we henceforth drop the index $k$, as well as 
$\rho$ in $\LL_n^{\rho}$ since we work with a fixed density $\rho>0$.
Finally, a simple computation gives an expression
for $\LL_n$
\be{eq2.1}
\LL_n(\v)= \LL_{(n)}(\v)
+\sum_{i\in \L_n}p^*_n(i,i) \gr(\v\circ\An^+_i -\v)+
\sum_{i\in \L_n}p_n(i,i) g_i(\v\circ\An^-_i -\v).
\ee
Note that by definition of $\LL_n$, the product of measures $\theta_{\gamma(\rho)}$
over sites of $\L_n$, which we denote either by $\nur^{\L_n}$ or simply by $\nur$,
is the invariant measure for $\LL_n$. Also, we have $\LL_n^*(\v)=
E_{\nur}[\LL^*(\v)|\FF_{\L_n}]$.  Finally, we omit the simple proof that
$\LL_n$ is a monotone irreducible process.

We denote by $E^n_{\eta}$ (resp. $E^{(n)}_{\eta}$)
the law of the Markov process generated
by $\LL_n$ (resp. $\LL_{(n)}$) with initial configuration $\eta$.
We denote by $\bar\LL_n:= 1_{\A^c} \LL_n$ 
(resp. $\bar\LL_{(n)}:= 1_{\A^c} \LL_{(n)}$) the process
killed on $\A$, and by $\bar S_t^n$ (resp. $\bar S_t^{(n)}$)
the associated semi-group. Note that if $\tau$ is the first
occurrence time of $\A$, then for $\v|_{\A}\equiv 0$ 
\[
\bar S_t^n(\v)(\eta)=E^n_{\eta}[\v(\eta_{t\wedge \tau})]=
E^n_{\eta}[\v(\eta_{t})1_{\tau>t}].
\]
\subsection{Approximating the killed process.}
\label{sec-killed}
The main uniqueness result is the following.
\bl{lem-killed} For any $\v\in \Lb$ with $\v|_{\A}\equiv 0$,
we have
\[
\forall t>0,\qquad
\lim_{n\to\infty}\int|\bar S_t^n(\v)-\bar S_t(\v)|d\nur=0.
\]
\el
\bpr
We first approximate $\{\tau>t\}$ by $\{\eta(t_i)\not\in \A,
\ i=0,\dots,k\}$ where $\{t_i\}$ is a regular subdivision of
$[0,t]$ of mesh $t/k$; we denote the latter event $\{\tau^{k}>t\}$.
Thus, we show in Step 1 that for each $k>0$, and $\v\in \Lb$ 
with $\v|_{\A}\equiv 0$
\be{eq-kill1}
\lim_{n\to\infty}\int|
E_{\eta}^n[1_{\{\tau^{k}>t\}}\v(\eta_t)]-
E_{\eta}^{(n)}[1_{\{\tau^{k}>t\}}\v(\eta_t)]|d\nur=0.
\ee
Since by Lemmas 2.3 and 2.6 of \cite{andjel}, 
we have the pointwise convergence
\be{eq-ind0}
E_{\eta}^{(n)}[1_{\{\tau^{k}>t\}}\v(\eta_t)]=
S_{t_1}^{(n)}\left(
1_{\A^c}S_{t_2}^{(n)}\left(1_{\A^c}\dots S_{t_{k+1}}^{(n)}
(\v)\right)\right)(\eta)
\stackrel{\scriptstyle{n \ra \infty}}{\longrightarrow}
E_{\eta}[1_{\{\tau^{k}>t\}}\v(\eta_t)],
\ee
we would conclude that
\be{eq-ind01}
\lim_{n\to\infty} \int\! |E_{\eta}^n[1_{\{\tau^{k}>t\}}\v(\eta_t)]-
E_{\eta}[1_{\{\tau^{k}>t\}}\v(\eta_t)]|d\nur=0.
\ee

In Step 2, we show that there is a constant $C$ independent of $n$
such that
\be{eq-kill2}
\int\!\!|E_{\eta}^n[1_{\{\tau^{k}>t\}}\v(\eta_t)]-
E_{\eta}^n[1_{\{\tau>t\}}\v(\eta_t)]|d\nur\le C\e.
\ee
Also, leaving $\A$ requires that all the particles
in {\it excess} escape $\S$ in a subinterval of length $t/k$. Thus,
the continuity properties of the infinite volume process
give 
\be{eq-ind02}
\lim_{k\to \infty}E_{\eta}[1_{\{\tau^{k}>t\}}\v(\eta_t)]=
E_{\eta}[1_{\{\tau>t\}}\v(\eta_t)],
\ee
and the proof is concluded once we combine (\ref{eq-ind01}), 
(\ref{eq-kill2}) and (\ref{eq-ind02}).

\noindent{\bf Step 1.}

First, we show by induction on $k$ (the number of
points in the subdivision of $[0,t]$) that 
there are two constants $C_k,C_k'$ such that for $\eta\not\in\A$
if we set $\d_n(i)=(p_n(i,i)+p^*_n(i,i))\a(i)$
\be{eq-ind1}
|E_{\eta}^n[1_{\{\tau^{k}>t\}}\v(\eta_t)]-
E_{\eta}^{(n)}[1_{\{\tau^{k}>t\}}\v(\eta_t)]|\le
C_k\sum_{i\in\L_n} \d_n(i)
\sum_{j=0}^{k-1}\int_0^t\!\! E^n_{\eta}[\eta_{s+s_j}(i)+C_k']ds,
\ee
where $s_0=0$ and $s_j=t_1+\dots+t_j$.

For $k=1$, we have $t_0=0$ and $t_1=t$, so that 
(\ref{eq-ind1}) reduces
to show that for $\eta\not\in \A$, there are $C_1,C_1'$ such that
\be{eq-ind2}
|S_t^n\v(\eta)-S_t^{(n)}\v(\eta)|\le 
C_1 \sum_{i\in\L_n} \d_n(i)
\int_0^t\!\!  E^n_{\eta}[\eta_{s}(i)+C_1']ds,
\ee
To obtain (\ref{eq-ind2}), we use an integration by parts
formula
\[
S_t^n\v(\eta)-S_t^{(n)}\v(\eta)=
\int_0^t S^n_{t-s}(\LL_n-\LL_{(n)})S_s^{(n)}\v(\eta) ds.
\]
Since $\v\in \Lb$, Lemma 2.2 of~\cite{andjel} implies that
for some constant $C$
\[
L(S_s^{(n)}\v)\le e^{Cs} L(\v).
\]
From (\ref{eq2.1}) it is enough to bound terms of the form
\be{eq-ind3}
|\An_i^{\pm}S_s^{(n)}\v(\eta)-S_s^{(n)}\v(\eta)|\le L(
S_s^{(n)}\v)\a(i)\le L(\v) e^{Cs} \a(i).
\ee
Thus, 
\[
|S_t^n\v(\eta)-S_t^{(n)}\v(\eta)|\le
L(\v) \sum_{i\in\L_n} \d_n(i)
\int_0^t\!\! \left(S^n_{s}(g_i)(\eta)+\gr\right)ds.
\]
(\ref{eq-ind2}) follows after recalling that $g_i(\eta)\le \Delta\eta(i)$.

The induction step from $k$ to $k+1$ follows with exactly the same
arguments. First, we recall (\ref{eq-ind0}) and write similarly
\[
E_{\eta}^{n}[1_{\{\tau^{k}>t\}}\v(\eta_t)]=
S_{t_1}^{n}\left( 1_{\A^c}S_{t_2}^{n}\left(1_{\A^c}\dots
S_{t_{k+1}}^{n}(\v)\right)\right).
\]
We call $\psi_2:=S_{t_2}^{(n)}(1_{\A^c}S_{t_3}^{(n)}
(1_{\A^c}\dots))$,
and recall that $\psi_2\in \Lb$ by Lemma 2.3 of \cite{andjel}. 
We now show that $1_{\A^c}\psi_2\in \Lb$. Indeed, for $\eta,\zeta\in\Omega$
\be{eq-kill4}
|\psi_2(\eta)1_{\eta\in\A^c}-\psi_2(\zeta)1_{\zeta\in\A^c}|\le
1_{\eta,\zeta\in\A^c}|\psi_2(\eta)-\psi_2(\zeta)|+
1_{B}(\eta,\zeta)|\psi_2|_{\infty}.
\ee
where we set $B:=\A\times\A^c\cup \A^c\times\A$.
Now, $(\eta,\zeta)\in B$ implies that $\sum_{\S}|\eta(i)-
\zeta(i)|\ge 1$. Thus, 
\be{eq-kill3}
1_B(\eta,\zeta)\le \sum_{\S}|\eta(i)- \zeta(i)|\le
\frac{\sum_{\S}|\eta(i)- \zeta(i)|\a(i)}{\inf_{\S} \a(i)}\le 
C||\eta-\zeta||.
\ee
Thus, combining (\ref{eq-kill3}) and (\ref{eq-kill4}) we 
obtain that $1_{\A^c}\psi\in \Lb$. Now,
\ba{eq-ind4}
E_{\eta}^n[1_{\{\tau^{k}>t\}}\v(\eta_t)]-
E_{\eta}^{(n)}[1_{\{\tau^{k}>t\}}\v(\eta_t)]&=&
\left(S_{t_1}^{n}(1_{\A^c}\psi_2)-S_{t_1}^{(n)}(1_{\A^c}\psi_2)\right)\cr
&&-S_{t_1}^{n}\left(
1_{\A^c}\left(\psi_2-S_{t_2}^{n}(1_{\A^c}S_{t_3}^{n}(
1_{\A^c}\dots))\right)
\right)
\ea
To the first term on the r.h.s we apply the estimates of the step $k=1$
of the induction. For the second term, the difference
$\psi_2-S_{t_2}^{n}(1_{\A^c} S_{t_3}^{n}(1_{\A^c}\dots))$ has $k$ subdivision
times, and we use our induction hypothesis to obtain (\ref{eq-ind1})
at order $k$; since $S_{t_1}^{n}$ is positive preserving, the
inequality is preserved after applying $S_{t_1}^{n}$ and we obtain
the desired (\ref{eq-ind1}) at order $k+1$. Now, to obtain (\ref{eq-ind0}),
note that 
\[
\sum_{i\in \L_n} \d_n(i)\le C\sum_{i\in \L_n\bs\L_{n-R}}\a(i)
\stackrel{\scriptstyle{n \ra \infty}}{\longrightarrow}
0\qquad(\text{since }\sum_{i\in \Z^d}\a(i)<\infty).
\]

\noindent{\bf Step 2.}
Let $\sigma_{\S}$ be the first time a particle inside $\S$ escapes $\S$,
and let $\theta_t$ be the time-translation by $t$. By the
strong Markov property, for $\eta\not\in\A$ and $\e=t/k$
\ba{eq-strong1}
|P_{\eta}^n(\tau>t)-P_{\eta}^n(\tau^{k}>t)|&\le&
P^n_{\eta}\left(\bigcup_{i\le k}
\{\tau\in [t_{i-1},t_i[,\sigma_{\S}\circ \theta_{\tau}<\e\}\right)\cr
&=&
\sum_{i=1}^k E_{\eta}^n[1_{\tau\in[t_{i-1},t_i[}P_{\eta_{\tau}}^n(\sigma_{\S}<\e)].
\ea
We need now a uniform estimate on $P_{\eta_{\tau}}^n(\sigma_{\S}<\e)\le C\e$.
By the hypotheses made on $\A$, we know that at time $\tau$, there is
a bounded number of particles in $\S$. For the zero range
process, it is routine to couple, from time $\tau$ onward,
the motion of the particle inside $\S$ (at time $\tau$) with a process 
containing only particles in $\S$ distributed as those of $\eta_{\tau}$. Now,
for this new process, at
any site, the rate of jump is bounded (uniformely in $\eta_{\tau}$,
since the number of particles is uniformely bounded), 
and the probability of having a jump before time $\e$
is smaller than $1-\exp(-\bar c\e)\le \bar c \e$. This concludes Step 2.
\epr
\subsection{Donsker-Varadhan functionals in $\L_n$}
\label{dv-n}
For $(\v,\mu)\in \D_n^{+}\times \M_n$, we define
\be{def-Gn}
\Gn(\v,\mu):=\int \frac{\LL_n \v}{\v}d\mu.
\ee
This is well defined since $\v>0$ on $\A^c$ which
contains the support of $\mu$. The functional $\Gn(\v,\mu)$
is useful if it has some regularity in $\mu$ and convexity in $\log(\v)$.
\bl{lem3.7} Assume $d\ge 3$. (i) For any $\v\in \D_n$, $\Gn(\v,.):\M_n\to \R$
is continuous. (ii) For any $\mu\in \M_n$, the map
$\tilde \Gn(.,\mu):=\Gn(\exp(.),\mu): \E_n\to \R$ is convex. 
\el
\bpr
Since $\LL_n(\v)/\v$ is not bounded, point (i) is not obvious. 
Let $\{\mu_k,k\in \N\}$ be in $\M_n$ converging weakly to $\mu$.
We show that for any $\v\in \D_n$, $\Gn(\v,\mu_k)$ converges to
$\Gn(\v,\mu)$ as $k$ tends to infinity. We recall the
notation $\nabla^i_j=T^i_j-1$,
\be{eq3.23}
\Gamma_n(\v,\mu_k)
:=\sum_{i,j\in \L_n}p(i,j)\int\!\!g_i \frac{\nabla^i_j \v}{\v} d\mu_k+
\sum_{i\in \L_n}\int\!\!( p^*_n(i,i)\gr \frac{\An^+_i \v-\v}{\v}+
p_n(i,i)g_i \frac{\An^-_i \v-\v}{\v}) d\mu_k.
\ee
Let $K_M$ be the compact set defined in (\ref{def-K}).
When integrating over $K_M$, the integrals on the r.h.s of 
(\ref{eq3.23}) pose no problem since the 
integrant over $K_M$ is bounded.
When integrating over $K_M^c$, first we recall that by
Lemma~\ref{lem3.3}, we have that $\v,1_{\A^c}/\v,g_i$
as well as $f_k:=d\mu_k/(d\nur)$
are in $L^p(\nur)$ for any $p\ge 1$. We then use H\"older's
inequality for $p=5$
\ba{eq3.26}
\int_{K_M^c}g_i\frac{T^i_j \v}{\v}d\mu_k
&&\le \int_{K_M^c}g_i\frac{\v\circ\An^-_i}{\v}f_kd\nur\cr
&&\le \left(\int g_i\v^p\circ\An^-_id\nur
\int\frac{1_{\A^c}}{\v^p}d\nur 
\int f_k^pd\nur\int g_i^{p-1}d\nur\nur(K_M^c)\right)^{1/p}\cr
&&\le \left(\gr\int \v^pd\nur
\int\frac{1_{\A^c}}{\v^p}d\nur \int f_k^pd\nur\int g^{p-1}d\nur\right)^{1/p}
\nur(K_M^c)^{1/p}\cr
&&\le C \nur(K_M^c)^{1/p}
\stackrel{\scriptstyle{M \ra \infty}}{\longrightarrow} 0.
\ea
The other terms of (\ref{eq3.23}) are dealt with in the same way.
To establish (ii), note first that by Lemma~\ref{lem3.2}, $\E_n$ is
convex. Then
\ba{eq3.20}
\Gamma_n(e^h,\mu)&=&\sum_{i,j\in \L_n}\!\!p(i,j)
\int g_i( e^{\nabla^i_j h}-1)d\mu\cr
&&+\sum_{i\in \L_n}p^*_n(i,i)\gr\int (e^{h\circ\An^+_i-h}-1)+
p_n(i,i)\int\!\!g_i(e^{h\circ\An^-_i-h}-1) d\mu.
\ea
The convexity follows from the convexity of the exponential.
\epr
\subsection{A variational formula for $\lnr$.}
\label{varia-ln}
\bl{lem3.4} For $d\ge 1$, there is $u_n\in \D_n$ and
$\lambda_n(\rho)>0$ such that 
\be{eq3.6}
1_{\A^c}\LL_n(u_n)+\lnr u_n=0.
\ee
Moreover $u_n$ is positive on $\A^c$.

Similary, when $d\ge 1$, there is $u_n^*\in \D_n^{*}$, positive on $\A^c$,
which satisfies $1_{\A^c}\LL_n^* u_n^*+\lnr u_n^*=0$, and
\be{eq3.7}
-\lnr=\lim_{t\to\infty} \frac{1}{t}\log(P^{n}_{\nur}(\tau>t)).
\ee
\el
\bpr
The proof follows the same lines as that of \cite{ac}
(see also \cite{fkmp}). This is
expected since $\LL_n$ is a monotone operator with the same features
as $\LL$. Thus, (\ref{eq3.7}) follows as simply as (\ref{eq0.6})
by a subadditivity argument. Now, for $\eta\in \Omega_n$, we denote
\be{eq3.8}
u_{t,n}(\eta)=\frac{P_{\eta}^{n}(\tau>t)}{P_{\nur}^{n}(\tau>t)}=
\frac{e^{t 1_{\A^c}\LL_n}(1_{\A^c})(\eta)}{P_{\nur}^{n}(\tau>t)},
\text{ and }
u^*_{t,n}(\eta)=\frac{e^{t 1_{\A^c}\LL^*_n}(1_{\A^c})(\eta)}
{P_{\nur}^{n}(\tau>t)}
\ee
and as in Step 1 of the proof of Lemma 2.6 of \cite{ac}, $u_{t,n}\in \D_n$
and $u^*_{t,n}\in\D_n^*$.
We focus now on $u_{t,n}$, though similar properties will hold
for $u^*_{t,n}$. First, by Lemma~\ref{lem-approx}, $\lnr\ge\lambda(\rho)>0$.
Thus, for any $k$, $\int_0^{\infty} P_{\nur}^{n}(\tau>t) t^k dt<\infty$, 
and as in Remark~\ref{rem-existence} we define
\[
dm_k(t)=\frac{P_{\nur}^{n}(\tau>t) t^k dt}{\int_0^{\infty}
P_{\nur}^{n}(\tau>t) t^k dt}\qquad\text{and}\qquad
\frac{d\Phi^{(k)}_n(\nur)}{d\nur}(\eta)=\int_0^{\infty}
u_{t,n}(\eta) dm_k(t).
\]
With identical arguments as in the proof of Theorem 2.4
of~\cite{ac}, the Cesaro weak-$L^2(\nur)$
limits of $\{\Phi_n^{(k)}(\nur),k\in \N\}$
are solutions of (\ref{eq3.6}). Now, it is clear that
$d\Phi^{(k)}_n(\nur)/d\nur\in\D_n$. Also, in the
weak-$L^2(\nur)$ topology $\D_n$ is compact by Remark~\ref{rem-comp},
and contain all
the Cesaro weak limits of $\{\Phi_n^{(k)}(\nur),k\in \N\}$. Thus, there is
a solution of (\ref{eq3.6}) in $\D_n$: we denote it by $u_n$.

We now show that $u_n>0$ on $\A^c$. By contradiction
assume that for $\eta\in\Omega_n\bs \A$, $u_n(\eta)=0$. Then 
(\ref{eq3.7}) implies that
$\LL_n(u_n)(\eta)=0$. This, in turn, implies that 
\begin{itemize}
\item[(i)] For all $i,j\in \L_n$ with $p(i,j)>0$, 
we have $u_n(T^i_j\eta)=0$.
\item[(ii)]For all $i\in \L_n$ with $p_n^*(i,i)>0$, we have
$u_n(\An_i^+\eta)=0$.
\item[(iii)] For all $i\in \L_n$ with $\eta(i)p_n(i,i)>0$, we have
$u_n(\An_i^-\eta)=0$.
\end{itemize}
To conclude that $u_n\equiv 0$ on $\A^c$, it is enough to note that by the
hypotheses $(\C\!\!-\!\!\FF)$ on $\A^c$, each $\eta\in\A^c$ can be transformed
into $O_{\L_n}$ by a succession of actions $\{\An_i^-\}$ with
$i\in\L_n$, and $\{T^i_j\}$ with $i,j\in\L_n$. The reverse
operation is made through a succession of $\{\An_i^+\}$ with
$i\in\L_n$, and $\{T^i_j\}$ with $i,j\in\L_n$.
\epr

We now establish the Donsker-Varadhan representation for
$\lnr$.
\bl{lem3.8} Assume $d\ge 3$. If $\A$ satisfies $(\C\!\!-\!\!\FF)$ of
Section~\ref{sec-pattern}, then $\lnr$ is given by
\be{eq3.19}
-\lnr=\sup_{\mu\in \M_n}\inf_{\v\in \D_n^{+}}
\int \frac{\LL_n \v}{\v} d\mu.
\ee
\el
\bpr
Let us call $\gamma_n$ the right hand side of (\ref{eq3.19}).
From Lemma~\ref{lem3.4}, there is $u_n\in \D_n^{+}$ such that
$\bar\LL_n u_n+\lnr u_n=0$. This implies that $\gamma_n\le -\lnr$.
We can use a classical minimax theorem \cite{fan}, since we have that
(i) for any fixed $\mu\in \M_n$, $h\mapsto \tilde\Gn(h,\mu)$ is convex
(by Lemma~\ref{lem3.7}) on the convex set $\E_n$ (by Lemma~\ref{lem3.2}),
(ii) for any fixed $h\in\E_n$, $\mu\mapsto \tilde\Gn(h,\mu)$ is continuous
(by Lemma~\ref{lem3.7}) on the compact set $\M_n$.
Thus,
\be{eq3.221}
\gamma_n=\inf_{\v\in \D_n^{+}}\sup_{\mu\in \M_n}
\int \frac{\LL_n \v}{\v} d\mu.
\ee
Now, for any $\v\in \D_n^{+}$, $0<\int\!\!\v u_n^*d\nur<\infty$, and
we can define
\[
d\mu^*=\frac{ \v u_n^*d\nur}{\int \v u_n^*d\nur}\in \M_n
\quad(\text{by Lemma~\ref{lem3.3}}). 
\]
Then, by duality
\[
\int\frac{\LL_n(\v)}{\v}d\mu^*=
\int \frac{\LL_n(\v)}{\v}\frac{\v u_n^*}{\int \v u_n^* d\nur} d\nur=
\int \frac{\v}{\int \v u_n^* d\nur}\LL_n^*(u_n^*)d\nur=-\lnr.
\]
By (\ref{eq3.221}), $\gamma_n\ge -\lnr$, and the proof is concluded.
\epr

In the following lemma, we establish
the uniqueness of the principal Dirichlet eigenfunction.
\bl{lem3.3bis} Assume $d\ge 3$. 
There is a unique non-negative eigenfunction $u_n\in \D_n$
of $1_{\A^c}\LL_n$ which satisfies $\int u_n d\nur=1$.
\el
\bpr
We know from Lemma~\ref{lem3.4} that there exists a 
positive eigenfunction $u_n$.
Assume that $\tilde u$ is a non-negative
Dirichlet eigenfunction with $\int \tilde u d\nur=1$ and
corresponding eigenvalue $\tilde \lambda$.
By the same argument as in the proof of Lemma~\ref{lem3.4}, we have that
$\tilde u$ is positive on $\A^c$.

First, we show that $\tilde \lambda=\lambda_n$.
Let $u_n^*$ be the dual eigenfunction given in Lemma~\ref{lem3.4}.
We multiply equality (\ref{eq3.6})
by $u_n^*$, integrate over $\nur$ and use duality
\be{eq3.12}
\int u_n^* \LL_n(\tilde u)d\nur=-\tilde \l
\int u_n^* \tilde u d\nur\Longrightarrow
(\lnr-\tilde \l)\int u_n^* \tilde ud\nur=0.
\ee
Now, since $u_n^*$ and $\tilde u$ are positive on $\A^c$ we conclude that
$\tilde \l =\lnr$.

Second, we show that $\tilde u=u_n$. Set
$h:=\log(u_n)$ and $\tilde h:=\log(\tilde u)$, on $\A^c$.
For any $\mu \in \M_n$ and any $\gamma\in ]0,1[$, by the convexity
of $\tilde \Gn$
\be{eq3.13}
\gamma\tilde \Gamma_n(h,\mu)+(1-\gamma) \tilde\Gamma_n(\tilde h,\mu)\ge
\tilde\Gamma_n(\gamma h+(1-\gamma) \tilde h,\mu).
\ee
Since $u_n$ and $\tilde u$ are solution of (\ref{eq3.6}), the left hand side
of (\ref{eq3.13}) is $-\lnr$. We define $h_{\gamma}=
\gamma h+(1-\gamma)\tilde h\in \E_n$ and we note that
$0<\int\exp(h_{\g}) u^*_n d\nur<\infty$. Now,
\[
d\mu_{\gamma}=\frac{e^{h_{\gamma}}u^*_nd\nur}{
\int\!\!e^{h_{\gamma}}u^*_nd\nur}\in \M_n,\qquad\text{and is such that}\quad
\tilde \Gamma_n(h_{\gamma},\mu_{\gamma})=
\Gamma_n(u_n^*,\mu_{\gamma})=-\lnr.
\]
Thus, we have equality in (\ref{eq3.13}) with $\mu_{\gamma}$.
Since $\mu_{\gamma}$ gives a positive weight to 
any $\eta\in \Omega_n\bs \A$, the following three conditions hold:
(i) for all $i,j\in \L_n$ with $g_i(\eta)p(i,j)>0$, we have $\nabla^i_j \tilde h=
\nabla^i_j h$; (ii)
for all $j\in \L_n$ with $p^*_n(j,j)>0$, we have
$(\An_j^+-\id)\tilde h=(\An_j^+-\id)h$;
(iii) for all $j\in \L_n$ with $g(\eta(j))p_n(j,j)>0$, we have
$(\An_j^--\id)\tilde h=(\An_j^--\id)h$.

Since $u_n$ is positive on $\A^c$, we form $f=\tilde u/u_n$, and rewrite
the conditions (i)-(iii) for $f$.
\begin{itemize}
\item[(i)] For all $\eta\in \A^c$ and $i,j\in \L_n$ with $\eta(i)p(i,j)>0$, 
we have $f(T^i_j\eta)=f(\eta)$.
\item[(ii)]For all $i\in \L_n$ with $p_n^*(i,i)>0$, and $\An_i^+\eta
\in \A^c$, we have $f(\An_i^+\eta)=f(\eta)$.
\item[(iii)] For all $\eta\in \A^c$ and 
$i\in \L_n$ with $\eta(i)p_n(i,i)>0$, we have
$f(\An_i^-\eta)=f(\eta)$.
\end{itemize}
As in the proof of Lemma~\ref{lem3.4}, we conclude that $\tilde u=u_n$.
\epr
\subsection{Approximating the principal eigenvalue.}
\label{approx-ln}
With an abuse of notations, we define for any finite
domain $U$, $\LL_U(\v)= E_{\nur}[\LL(\v)|\FF_U]$. We mean
by $\LL_U$ an expression like (\ref{eq2.1}) where $U$ replaces
$\L_n$: thus, a zero-range process on $U$ with creations and annihilations
on the boundaries of $U$. We denote by $S^U_t$ the semi-group
associated with $\LL_U$ and by $P^U_{\nu}$ the corresponding Markov
process with initial measure $\nu$.
We denote by $\bar S^U_t$ the semi-group killed on $\A$. 

We first state an obvious corollary of Lemma~\ref{lem-killed}
applied to $\v=1_{\A^c}$.
\bc{cor1}
When the pattern satisfies $(\C\!\!-\!\!\FF)$, we have
\[
\lim_{n\to\infty} P^n_{\nur}(\tau>t)=P_{\nur}(\tau>t).
\]
\ec

\noindent{\bf Proof of Lemma~\ref{lem-approx}}
We divide the proof in two steps.

\noindent{\bf Step 1.} We show that $n\mapsto P^{n}_{\nur}(\tau>t)$
is increasing.

Let $U$ be a finite subset, $i\not\in U$, and set $\tilde U=U\cup \{i\}$.
Thus, it is enough to show that $\int (\bar S^{\tilde U}_t 1_{\A^c}-
\bar S^{U}_t 1_{\A^c})d\nur\ge 0$. Step 1 follows then by
induction. Note that for $\v$ $\FF_U$-measurable and $j\in U$, we have
$\v\circ T^i_j=\v\circ\An^+_j$, $\v\circ T^j_i=\v\circ \An^-_j$,
$\v\circ \An_i^+=\v$ and $\v\circ\An_i^-=\v$ so that 
\ba{eq2.4}
(\bar \LL_{\tilde U}-\bar \LL_{U})\v &=& 1_{\A^c}
\sum_{j\in U}\left( p(j,i) g_j(\v\circ T^j_i-\v)+
p(i,j) g_i(\v\circ T^i_j-\v)\right)\cr
&&-1_{\A^c} \sum_{j\in U}\left( p(j,i) g_j(\v\circ \An^-_j-\v)+
p(i,j) \gr(\v\circ \An^+_j-\v)\right)\cr
&=&1_{\A^c} \sum_{j\in U}
p(i,j)(g_i- \gr)\left(\v\circ \An^+_j-\v\right).
\ea
Now, we set $\v_s:=\bar S^U_s(1_{\A^c})$ and $\psi_s:
=(\bar S^{\tilde U}_s)^* (1_{\A^c})$, 
and we use an integration by parts formula
\ba{eq2.5}
\int\!\!\bar S^{\tilde U}_t(1_{\A^c})d\nur-
\int\!\!\bar S^{U}_t(1_{\A^c})d\nur &=&
\int\!\!\!\int_0^t
\bar S^{\tilde U}_{t-s}(\bar \LL_{\tilde U}-\bar \LL_U)
\bar S^U_s(1_{\A^c})dsd\nur\cr
&=&\int\!\!\!\int_0^t (\bar \LL_{\tilde U}-\bar \LL_U)(\v_s)
\psi_{t-s}dsd\nur.
\ea
Thus, by (\ref{eq2.4})
\ba{eq2.6}
P^{\tilde U}_{\nur}(\tau>t)&-&P^{U}_{\nur}(\tau>t)=
\!\sum_{j\in U} p(i,j)\int\!\!\!\int_0^t\!\!(\An^+_j\v_s-\v_s)
(g_i-\gr)\psi_{t-s}  dsd\nur^{\tilde U}\cr
&=&\!\sum_{j\in U} p(i,j)\int\!\!\!\int_0^t\!\!(\An^+_j\v_s-\v_s)
\int\!\!(g_i-\gr)\psi_{t-s}d\nur^{\{i\}} dsd\nur^{U}.
\ea
Note that for any $s$,
$\eta\mapsto \psi_s(\eta)$ and $\eta\mapsto \v_s(\eta)$ is
decreasing positive, whereas $\eta\mapsto g_i(\eta)$ is increasing
and $\int g_id\nur=\gr$. Thus, by FKG inequality
\be{eq2.7}
\int (g_i-\gr)\psi_{t-s}d\nur{\{i\}}\le
\int (g_i-\gr)d\nur^{\{i\}}
\int \psi_{t-s} d\nur^{\{i\}}=0.
\ee
Thus, as $\v_s\circ\An^+_j-\v_s\le 0$, the first step concludes.
We call $\lambda_{\infty}(\rho)$ the limit of $\lambda_n(\rho)$.

\noindent{\bf Step 2.} We show the following Lemma which allows
us the conclude the proof of Lemma~\ref{lem-approx} readily.
\bl{lem6} Any subsequence of $\{u_n\}$ has a further subsequence
converging, in weak-$L^2(\nur)$,
to a solution $u$ of (\ref{eq0.7}(i)), and $u\in \D_{\rho}$.
Moreover, $\lambda_{\infty}(\rho)=\lr$.
\el
\bpr
For notational convenience, we write the proof for $\{u_n^*\}$.
Recall that $\D^*_{\rho}\cap\{\v:\ \int\!\!\v d\nur=1\}$ is compact
in weak-$L^2(\nur)$ by Remark~\ref{rem-comp}. Let $u^*\in\D^*_{\rho}$
be a limit point of $\{u_n^*\}$ along a subsequence which for simplicity
we still call $\{u_n^*\}$.
For any $\v\in \Lb$, and any integer $n$
\be{eq6.12}
\int \bar S_t^n(\v) u_n^*d\nur=e^{-\lnr t}\int \v u_n^* d\nur.
\ee
Then,
\ba{eq6.13}
|\int\!\! \bar S_t^n(\v) u_n^*d\nur&-&\int\!\! \bar S_t(\v) u^*d\nur|=
|\int\left( \bar S_t^n(\v)-\bar S_t(\v)\right) u_n^* d\nur|+\!\!
|\int \!\!\bar S_t\v(u_n^*-u^*) d\nur|\cr
&\le& \sup_n ||u_n^*||_{\nur} ||\bar S_t^n(\v)-\bar S_t(\v)||_{\nur}+\!\!
|\int \!\!\bar S_t\v(u_n^*-u^*) d\nur|.
\ea
The $L^2(\nur)$ convergence of $\bar S_t^n(\v)-\bar S_t(\v)$ is equivalent
to an $L^1(\nur)$ convergence, since $\v$ is bounded and $\bar S_t,
\bar S_t^n$ are contractions (in $L^{\infty}$). By recalling 
Lemma~\ref{lem-killed} and Step 1, and taking the limit $n$ to infinity, 
\be{eq6.14}
\int \bar S_t(\v) u^*d\nur=e^{-\lambda_{\infty}(\rho) t}\int \v u^* d\nur.
\ee
Now, since $\Lb$ is a dense set in $L^2(\nur)$, this implies that
$u^*\in \D(\bar\LL^*, L^2(\nur))$, and that (\ref{eq6.14}) holds
for any $\v\in L^2(\nur)$. Take $\v=u\in \D_{\rho}\subset
L^2(\nur)$ solution of (\ref{eq0.7}(i)), and use that
\[
\bar S_t(u)=e^{-\lr t} u,\quad\nur-\text{a.s.}
\Longrightarrow (e^{-\lambda_{\infty}(\rho) t}-e^{-\lr t})
\int\!\! u u^* d\nur=0.
\]
Now, since $u$ and $u^*$ are decreasing, and in $L^2(\nur)$, we have
\[
\infty>||u||_{\nur}||u^*||_{\nur}\ge \int\!\! u  u^* d\nur
\stackrel{\scriptstyle{\text{FKG}}}{\ge}
\int\!\! ud\nur\int\!\! u^* d\nur=1.
\]
Thus, $\lambda_{\infty}(\rho)=\lr $, and $u^*$ satisfies (\ref{eq0.7}(ii)).
\epr
\section{Donsker-Varadhan functionals on $\Z^d$}
\label{sec-dv}
The main problem arises since $\LL(\v)$ does not make sense
as a pointwise convergent series when $\v\in \D_{\rho}$. Indeed,
even if $\v$ were bounded, the naive bound $|\nabla^i_j\v|\le |\v|_{\infty}
(\e_i+\e_j)$ would fail since $\sum_k \e_k=\infty$.
Thus, we show in this section how to
obtain $\Gi(\v,\mu)$ as the limit of the Cauchy sequence
$\{\int\!\!\LL_n^{\rho}(\v)/\v d\mu,\ n\in \N\}$ taking advantage
of the gradient bounds on $\v$ and $d\mu/d\nur$ by an integration
by parts formula.
\subsection{Technical prerequisites}
We first define a family of functionals, $\{\Gn,\ n\in \N\}$,  
on $\D_{\rho}^+\times \M_{\rho}$,
whose limit when $n$ tends to infinity is shown to exist.
\bl{lem4.1} Assume $d\ge 3$. For $\v\in \D_{\rho}^+$ and $\mu\in \M_{\rho}$, 
and any integer $n$, the functional 
$\Gn(\v,\mu):=\int \LL_n (\v)/\v d\mu$ is well defined.
If we call $\tilde \Gn(h,\mu):=\Gn(\exp(h),\mu)$, then
for any $\mu$, 
the map $h\mapsto \tilde \Gn(h,\mu)$ is convex on the
convex set $\E_{\rho}$.
\el
\bpr
The formal full expression of $\Gn(\v,\mu)$ is
\be{eq4.0}
\Gn(\v,\mu)=
\sum_{i,j\in \L_n} p(i,j)\int g_i\frac{\nabla^i_j \v}{\v}d\mu
+\sum_{i\in \L_n}\left(\gr p_n^*(i,i)\int \frac{\nabla_i^+\v}{\v}d\mu
+p_n(i,i)\int g_i \frac{\nabla_i^-\v}{\v}d\mu\right).
\ee
Note that as $\v\in \D_{\rho}^+$, $T^i_j\v\le \An_i^-\v$. Thus,
(\ref{eq4.0}) is defined if we bound $\int g_i \An_i^{-}(\v)/\v d\mu$
for each site $i\in \L_n$. This is done as in (\ref{eq3.26}).

From (\ref{eq4.0}), an expression for $\tilde \Gn(h,\mu)$
is as follows
\ba{eq4.0tilde}
\tilde \Gn(h,\mu)&=&
\sum_{i,j\in \L_n} p(i,j)\int g_i\left(e^{\nabla^i_j h}-1\right) d\mu\cr
&&\qquad+\sum_{i\in \L_n}\left(\gr p_n^*(i,i)\int (e^{\nabla_i^+h}-1)d\mu
+p_n(i,i)\int g_i (e^{\nabla_i^-h}-1)d\mu\right).
\ea
The convexity of $h\mapsto \tilde \Gn(h,\mu)$ follows from
the convexity of the exponential.
\epr

We now express $\tilde \Gn(h,\mu)$ in terms of gradients of
$h$ and $\mu$. 

\bl{lem.gradient} For $h\in \E_{\rho}$ and
$\mu\in \M_{\rho}$, we have with $f:=d\mu/d\nur$
\be{eq4.7}
\tilde \Gn(h,\mu)= \sum_{i,j\in \L_n}
\gr p_n(i,j)
\left(
\int\!\left(e^{\Delta_i^j h}-1\right) \nabla^+_i f d\nur+ 
\int\!\left(e^{\Delta_i^j h}-1-\Delta_i^j h\right)d\mu
\right)+R_n(h,\mu),
\ee
with,
\be{eq4.1}
\lim_{n\to\infty} \sup_{h\in \E_{\rho}}\sup_{\mu\in \M_{\rho}} |R_n(h,\mu)|=0.
\ee
Note also that
\be{eq4.7bis}
\int \LL_n(\v)d\mu=
\sum_{i,j\in \L_n}\gr p_n(i,j)
\int\Delta_i^j \v \nabla^+_i f d\nur-\sum_{i\in \L_n}\gr p_n(i,i)
\int\nabla_i^+ \v \nabla^+_i f d\nur.
\ee
\el
\bpr
First, we apply the integration by parts formula (\ref{eq-byparts})
to (\ref{eq4.0tilde}): 
\ba{eq4.6}
\tilde \Gn(h,\mu)&=&\sum_{\substack{i,j\in \L_n\\i\not= j}}
p(i,j)\int\! g_i\An_i^-\left(e^{\Delta_i^j h}-1\right) d\mu\cr
&&\qquad+\sum_{i\in \L_n}\left(p_n^*(i,i) \gr\int\!(e^{\nabla^+_i h}-1) d\mu+
p_n(i,i)\int\!g_i\An_i^-(e^{-\nabla^+_i h}-1) d\mu\right)\cr
&=&\gr\!\!\sum_{i,j\in \L_n} p_n(i,j)
\left(\int\! \left(e^{\Delta_i^j h}-1\right) \nabla^+_i f d\nur
+\int\left(e^{\Delta_i^j h}-1-\Delta_i^j h\right)d\mu\right)\cr
&&\qquad\qquad+R_n(h,\mu)+N(h,\mu),
\ea
with 
\ba{eq4.8}
R_n(h,\mu)&:=&\sum_{i\in \L_n}\gr p_n(i,i)
\int\!(e^{-\nabla^+_i h}-1)\nabla^+_i (f) d\nur\cr
&&+\gr\sum_{i\in \L_n}\!\! p^*_n (i,i) \int\!(e^{\nabla^+_i h}-1
-\nabla^+_i h) d\mu\cr
&&+\gr\sum_{i\in \L_n}\!\! p_n (i,i) \int\!(e^{-\nabla^+_i h}-1
+\nabla^+_i h) d\mu,
\ea
and,
\be{eq4.9}
N(h,\mu):=\sum_{i,j\in \L_n}\gr p(i,j)\int \Delta_i^j hd\mu+
\sum_{i\in \L_n}\gr(p^*_n(i,i)-p_n(i,i))\int(\nabla_i^+h)d\mu.
\ee
To show that $N(h,\mu)$ vanishes, first write
\ba{eq-mistaken}
\sum_{i,j\in \L_n}p(i,j)\int(\nabla_j^+h)d\mu&=&
\sum_{j\in\L_n}\left(
\sum_{i\in\L_n}p(i,j)\right)\int (\nabla_j^+h) d\mu\cr
&=&\sum_{j\in\L_n}\left(1-
\sum_{i\not\in\L_n}p(i,j)\right)\int (\nabla_j^+h) d\mu\cr
&=&\sum_{j\in\L_n}(1-p_n^*(j,j))\int (\nabla_j^+h) d\mu,
\ea
and similarly,
\[
\sum_{i,j\in \L_n}p_n(i,j)\int(\nabla_i^+h)d\mu=
\sum_{i\in\L_n}(1-p_n(i,i))\int (\nabla_i^+h) d\mu.
\]
It is thus clear that $N(h,\mu)=0$.

We now show that $R_n(h,\mu)$ defined in (\ref{eq4.8}) is negligeable.
Note that for $i\not\in\S$, $\e_i< 1$. Also, for $n$ large
enough, if we define $\p^R\L_n:=\L_n\bs\L_{n-R}$, then
$\p^R \L_n\cap \S=\emptyset$. Also, by (\ref{def-p}) (iii),
$p_n(i,i)=0$ when $i\not\in \p^R\L_n$. 
Thus, there is a constant $c_0>0$ such that for $i\in \p^R\L_n$
$|\nabla_i^+ h|\le -\log(1-\e_i)\le c_0 \e_i$, $\nur$-a.s., and
$|\nabla_i^+ f|\le (\e_i+\e^*_i)f$, $\nur$-a.s.\ .Thus, there is
a constant $c_1>0$ such that
\be{eq4.2bis}
p_n(i,i)|\int \!\left( e^{-\nabla^+_i h}-1\right)\nabla^+_i f d\nur|
\le \int (\exp(c_0\e_i)-1)(\e_i+\e_i^*) d\mu \le c_1 \e_i(\e_i+\e_i^*),
\ee
and by expanding to second order in $\nabla_i^+h$
\be{eq4.2}
p_n^*(i,i)|\int(e^{\nabla_i^+h}-1-\nabla_i^+h)d\mu|\le c_1 \e_i^2,
\quad\text{and}\quad
p_n(i,i)|\int(e^{-\nabla_i^+h}-1+\nabla_i^+h)d\mu|\le c_1 \e_i^2,
\ee
Combining (\ref{eq4.2bis}) and (\ref{eq4.2}), and summing over
$i\in \p^R\L_n$, we obtain the desired asymptotics (\ref{eq4.1}),
since for dimension $d\ge 3$, $\sum \e_i^2<\infty$.

We obtain (\ref{eq4.7bis}) from (\ref{eq4.7})
by setting $h=\e \v$ and expanding
$\tilde \Gn(h,\mu)$ to first order in $\e$.
\epr

We are now ready for the key technical lemma of this section.
\bp{prop4.1} For $(\v,\mu)\in \D_{\rho}^+\times \M_{\rho}$,
$\{\Gn(\v,\mu),n\in \N\}$ is a Cauchy sequence whose
limit we denote by $\Gi(\v,\mu)$. We have
the following properties.
\begin{itemize}
\item[(i)] For $h\in \E_{\rho}$,
$h\mapsto \tilde\Gi(h,\mu):=
\Gi(e^h,\mu)$ is convex.
\item[(ii)] The Cauchy sequence is uniform in the following sense
\be{eq4.4bis}
\lim_{n\to\infty} \sup_{\v\in \D_{\rho}^+}\sup_{\mu\in \M_{\rho}}
|\Gn(\v,\mu)-\Gi(\v,\mu)|=0.
\ee
\item[(iii)] For any integer $n$, and any $\mu\in \M_{\rho}$ 
we denote by $\mu_n$ the measure of $\M_n$ of density 
$f_n:=E_{\nur}[d\mu/d\nur|\FF_{\L_n}]$. Then,
\be{eq4.4}
\lim_{n\to\infty} \sup_{\mu\in \M_{\rho}}\sup_{\v_n\in \D_n}
|\Gi(\v_n,\mu_n)-\Gi(\v_n,\mu)|=0.
\ee
\item[(iv)] For $\v_n\in \D_n^+$ and $\mu_n\in \M_n$, we have
$\Gi(\v_n,\mu_n)=\Gamma_{n}(\v_n,\mu_n).$
\end{itemize}
\ep
\bpr
\noindent{Step 1}:
We show that $\{\Gn(\v,\mu),n\in \N\}$ is a Cauchy sequence and
(\ref{eq4.4bis}) holds.

By using the expression (\ref{eq4.7}) of
Lemma~\ref{lem.gradient}, we have for
$m>n$
\ba{eq4.10}
\tilde \Gm(h,\mu)-\tilde \Gn(h,\mu)&=&
\sum_{\substack{(i,j)\in \L_m^2\bs \L_n^2\\i\not= j}}\gr
p(i,j)\left(\int\!(e^{\Delta^j_i h}-1) \nabla^+_i f d\nur+ 
\int\!(e^{\Delta^j_i h}-1-\Delta^j_i h) d\mu\right)\cr
&&\qquad+R_m(h,\mu)-R_n(h,\mu).
\ea
Since $p(i,j)=0$ when $|i-j|>R$, we can assume $n$ and $m$ so large
that if $(i,j)\in \L_m^2\bs \L_n^2$ with $p(i,j)>0$, then $i,j\not\in\S$.
Thus, there is a positive constant $c_0$ such that $\nur$-a.s.
\be{eq4.11} 
\forall (i,j)\in \L_m^2\bs \L_n^2\text{ with }p(i,j)>0,\quad
|\nabla^+_i h|\le -\log(1-\e_i)\le c_0 \e_i,\quad\text{and}\quad
|\nabla^+_i f|\le (\e_i+\e_i^*)f.
\ee
Also, there is a positive constant $c_1$ such that
\ba{eq4.13}
p(i,j)|\int(e^{\Delta^j_i h}-1) \nabla^+_if d\nur|&&\le
\int(e^{c_0(\e_i+\e_j)}-1)(\e_i+\e_i^*)f d\nur\cr
&&\le c_1(\e_i+\e_j)(\e_i+\e_i^*).
\ea
Now, recalling that for $i\not\in \S$, $\sum_j p(i,j)\e_j=\e_i$, 
and $\sum_j p(i,j)=1$, we have
\be{eq4.12}
\sum_{(i,j)\in \L_m^2\bs \L_n^2} c_1p(i,j) (\e_i+\e_j)(\e_i+\e_i^*)
\le 2c_1\sum_{i\in \L_n^c\cup \p^R \L_n}
\e_i(\e_i+\e_i^*) 
\stackrel{\scriptstyle{n \ra \infty}}{\longrightarrow}0,
\ee
since $\sum_i \e_i^2=\sum_i (\e_i^*)^2<\infty$ when $d\ge 3$.
Similarly, the second integral in (\ref{eq4.10}) will go to
0, after we perform a second order expansion and use (\ref{eq4.11}).
Now, from Lemma~\ref{lem.gradient}, $|R_m(h,\mu)-R_n(h,\mu)|$ converges
to 0 uniformely in $\E_{\rho}$ and $\M_{\rho}$. 

\noindent{Step 2}: The limit $h\mapsto \tilde\Gi(h,\mu)$ is
convex, since it is a pointwise limit of convex functions.

\noindent{Step 3}: We prove (\ref{eq4.4}). 

Let $\v_n$ be in $\D_n$, and
set $h_n=\log(\v_n)$. Note that for $i\not\in \L_n$, $\nabla_i^+ h_n=0$.
Also, for any function $\psi$, $\An_i^+E_{\nur}[\psi|\FF_{\L_n}]=
E_{\nur}[\An_i^+ \psi|\FF_{\L_n}]$. Thus, for $m>R+n$
\ba{eq4.15}
\tilde\Gm(h_n,\mu)&=&\tilde\Gn(h_n,\mu)+
\sum_{ \substack{i\in \L_m\bs \L_n\\j\in \L_n}}\!\!
\gr p(i,j)\int\left[ (e^{\nabla_j^+h_n}-1)\nabla_i^+ f+
(e^{\nabla_j^+h_n}-\nabla_j^+h_n-1)f\right] d\nur\cr
&&+ \sum_{ \substack{j\in \L_m\bs \L_n\\i\in \L_n}}\!\!
\gr p(i,j)\int\left[ (e^{-\nabla_i^+h_n}-1)\nabla_i^+ f+
(e^{-\nabla_i^+h_n}+\nabla_i^+h_n-1)f\right] d\nur\cr
&&+R_m(h_n,\mu)-R_n(h_n,\mu).
\ea
By observing that $\tilde\Gn(h_n,\mu)=\tilde\Gn(h_n,\mu_n)$,
and that $R_m(h_n,\mu)=0$ for $m>n+R$, we have
\ba{eq4.16}
\tilde\Gm(h_n,\mu)&-&\tilde\Gn(h_n,\mu_n)=
\sum_{ \substack{i\in \L_m\bs \L_n\\j\in\p^R \L_n}}
\gr p(i,j)\int(e^{\nabla_j^+h_n}-1)\nabla_i^+ fd\nur\cr
&&\sum_{j\in \p^R\L_n}\!\! \gr p_n(j,j)\left(
\int\!\!(e^{-\nabla_j^+h_n}-1)\nabla_j^+(f)d\nur
+\int\!\!(e^{-\nabla_j^+h_n}+\nabla_j^+h_n-1)f d\nur\right)\cr
&&\sum_{j\in \p^R\L_n}\!\! \gr p_n^*(j,j)
\int\!\!(e^{\nabla_j^+h_n}-\nabla_j^+h_n-1)f d\nur 
-R_n(h_n,\mu) .
\ea
Now, using again that for $j\in\p^R \L_n$,
$|\nabla_j^+ h_n|\le c_0 \e_j$, and for $i\in \p^R \L_n\cup\L_n^c$,
$|\nabla_i^+ f|\le (\e_i+\e_i^*) f$, we have a constant $C_1$ such that
\ba{eq4.17}
|\tilde\Gm(h_n,\mu)-\tilde\Gn(h_n,\mu_n)|&=&
|R_n(h_n,\mu)|+C_1\!\!\! 
\sum_{ \substack{i\in \L_m\bs \L_n\\j\in \L_n}}
\!\!\!p(i,j)\e_j(\e_i+\e_i^*)\cr
&&\qquad+2C_1\!\!\!
\sum_{j\in \p^R \L_n}\e_j^2+C_1\!\!\!\sum_{j\in \p^R \L_n}
\e_j(\e_j+\e^*_j)\cr
&\le& |R_n(h_n,\mu)|+2C_1\!\!\!
\sum_{i\not\in \L_n}\!\!\!\e_i^2+(\e_i^*)^2+
C_1\!\!\!  \sum_{j\in \p^R \L_n}(4\e_j^2+(\e^*_j)^2)
\ea
Equation (\ref{eq4.4}) follows after we take the limit
$m$ to infinity in (\ref{eq4.17}) and use (\ref{eq4.1}) of 
Lemma~\ref{lem.gradient}. 

\noindent{Step 4}: We show that
$\Gi(\v_n,\mu_n)=\Gn(\v_n,\mu_n)$. Indeed,
for $m>R+n$, $\LL_{(m)}(\v_n)=\LL(\v_n)$ so that
\be{eq4.18}
\Gm(\v_n,\mu_n)=
\int E_{\nur}\left[\frac{\LL_{(m)}(\v_n)}{\v_n}f_n\Big| \FF_{\L_m}\right]d\nur=
\int \frac{E_{\nur}[\LL_{(m)}(\v_n)|\FF_{\L_n}]}{\v_n}f_nd\nur=
\Gn(\v_n,\mu_n).
\ee
\epr

Now, a minimax theorem for $\Gi$ will be a corollary of 
Lemma~\ref{lem-compact}.
\bp{prop4.2}
A minimax theorem holds for $\Gi$. In other words,
\be{min-max}
\sup_{\mu\in \M_{\rho}} \inf_{\v\in \D_{\rho}}
\Gi(\v,\mu)=
\inf_{\v\in \D_{\rho}}\sup_{\mu\in \M_{\rho}} \Gi(\v,\mu).
\ee
\ep
\bpr
We need to check that for any $\v \in \D_{\rho}$, the map 
$\mu\mapsto \Gi(\v,\mu)$ on $\M_{\rho}$ is
continuous on the compact space $\M_{\rho}$. 
Let $\{\mu_k,k\in \N\}$ be in $\M_{\rho}$, converging
weakly to $\mu\in\M_{\rho}$. By Lemma~\ref{lem5.1}, all
densities $f_k=d\mu_k/d\nur$ are uniformely bounded in $L^2(\nur)$.
Thus $f_k$ converges in weak-$L^2(\nur)$ to $d\mu/d\nur$.
Now, for $\v\in \D^+_{\rho}$, as in (\ref{eq3.26}),
$g_i (\v\circ T^i_j)/\v\in L^2(\nur)$, so that for
$i,j\in \L_n$
\[
\int\!\! g_i \frac{\v\circ T^i_j}{\v}d\mu_k
\stackrel{\scriptstyle{k \ra \infty}}{\longrightarrow} 
\int\!\! g_i \frac{\v\circ T^i_j}{\v}d\mu.
\]
Thus, $\Gn(\v,\mu_k)\to\Gn(\v,\mu)$ as $k\to\infty$. Now, the
uniform Cauchy property (\ref{eq4.4bis}) implies that
$\Gi(\v,\mu_k)\to\Gi(\v,\mu)$ as $k\to\infty$.
\epr

\subsection{Proof of Theorem~\ref{the2}}
\label{sec-the2}
If $u_n$ is the principal normalized eigenfunction of
$\LL_n$, then for any $n$ and any $\mu_n$, we have by 
Proposition~\ref{prop4.1} (iv)
\be{eq4.19}
\Gi(u_n,\mu_n)=-\lambda_n(\rho).
\ee
Now, by (\ref{eq4.4}) of Proposition~\ref{prop4.1}, for any $\e>0$,
there is $n_0$ such that for any $n\ge n_0$
\be{eq4.20}
\sup_{\mu\in \M_{\rho}} |\Gi(u_n,\mu)-
\Gi(u_n,\mu_n)|\le \e.
\ee
Thus, for any $\mu\in \M_{\rho}$ and $n\ge n_0$
\be{eq4.21}
\Gi(u_n,\mu)\le -\lnr +\e\Longrightarrow
\inf_{\v\in \D_{\rho}}\Gi(\v,\mu)\le -\lnr +\e\qquad
(\text{since }\D_n\subset\D_{\rho}).
\ee
Recalling Lemma~\ref{lem-approx}, and taking the limit $n\to\infty$,
we obtain
\be{eq4.22}
\sup_{\mu\in \M_{\rho}}\inf_{\v\in \D_{\rho}}\Gi(\v,\mu)\le -
\lambda(\rho).
\ee
Conversely, if $h_n=E_{\nur}[h|\FF_n]$ for $h\in \E_{\rho}$, we show
by a convexity argument that
\be{eq4.23}
\forall \mu_n\qquad\tilde\Gi(h,\mu_n)\ge \tilde\Gi(h_n,\mu_n)-\e_n
\quad\text{with}\quad \lim_{n\to\infty} \e_n=0.
\ee
Indeed, take $m>n$ and in expression (\ref{eq4.7}) break down the
gradient $\nabla_i^+f_n$ so as to obtain
\be{eq4.24}
\tilde\Gm(h,\mu_n)= \sum_{i,j\in \L_m}
\gr p(i,j)\int \left[(e^{\Delta^j_ih}-1)\An_i^+ f_n-
(\Delta^j_ih)f_n\right] d\nur+R_{m}(h,\mu_n).
\ee
We further divide the sum over $\L_m$ into
\be{eq4.241}
\tilde\Gm(h,\mu_n)= \sum_{i,j\in \L_n}
\gr p(i,j)\int \left[(e^{\Delta^j_ih}-1)\An_i^+ f_n-
\Delta^j_ih)f_n\right] d\nur+Q_{n,m}(h,\mu_n),
\ee
where $Q_{n,m}(h,\mu_n)$ contains the sum over $(i,j)\in
\L_m^2\bs \L_n^2$. With similar estimates as those showing
that $R_n(h,\mu)$ goes to 0 when $n$ tends to infinity uniformely in
$h$ and $\mu$, in
the proof of Proposition~\ref{prop4.1}, $Q_{n,m}(h,\mu)$
goes to 0 as $n$ and $m$ tend to infinity.
Using that for any function $\psi$ and $i\in \L_n$,
$\An_i^+ E_{\nur}[\psi|\FF_n]=E_{\nur}[\An_i^+\psi|\FF_n]$, we have
by Jensen's inequality for the conditional expectation
\ba{eq4.25}
\tilde\Gm(h,\mu_n)&=& 
\sum_{i,j\in \L_n} \gr p(i,j)\int 
E_{\nur}[e^{\Delta^j_ih}-1|\FF_n]\An_i^+ f_n-
E_{\nur}[\Delta^j_ih|\FF_n]f_n d\nur
+Q_{n,m}(h,\mu_n)\cr
&\ge& \sum_{i,j\in \L_n} \gr p(i,j)\int 
(e^{\Delta^j_ih_n}-1)\An_i^+ f_n-
\Delta^j_i(h_n)f_n d\nur+Q_{n,m}(h,\mu_n)\cr
&\ge&
\tilde\Gn(h_n,\mu_n)-Q_n(h_n,\mu_n)+Q_{n,m}(h,\mu_n)\cr
&=& \tilde\Gi(h_n,\mu_n)-Q_n(h_n,\mu_n)+Q_{n,m}(h,\mu_n).
\ea
Thus, by taking the limit as $m$
tends to infinity, we obtain (\ref{eq4.23}) with
\be{eq4.26}
\e_n:=\lim_{m\to\infty} \sup_{h,\mu} (|R_{n,m}(h,\mu)|+|Q_n(h_n,\mu_n)|)
\stackrel{\scriptstyle{n \ra \infty}}{\longrightarrow} 0.
\ee
Now, for any $h\in \E_{\rho}$, since $\infty>\int \exp(h_n)u_n^*d\nur>0$, we
can define
\[
\frac{d\mu_n^*}{d\nur}=\frac{e^{h_n}u_n^*}{\int e^{h_n}u_n^*d\nur}.
\]
Thus, by duality 
$\tilde\Gamma_{n}(h_n,\mu_n^*)=\Gn^*(u_n^*,\mu_n^*)=-\lnr$.
and,
\be{eq4.27}
\sup_{\mu\in \M_{\rho}}
\tilde\Gi(h,\mu)\ge
\tilde\Gi(h,\mu_n^*)\ge 
\tilde\Gi(h_n,\mu_n^*)-\e_n= -\lnr-\e_n.
\ee
Thus, by taking the limit $n$ to infinity, 
and using Lemma~\ref{lem-approx}, we obtain
\be{eq4.28}
\sup_{\mu\in \M_{\rho}}\Gi(e^{h},\mu)\ge
-\lambda(\rho)\Longrightarrow
\inf_{\v\in \D^+_{\rho}} \sup_{\mu\in \M_{\rho}}\Gi(\v,\mu)\ge
-\lambda(\rho).
\ee
Now, since by Proposition~\ref{prop4.2}, the
minimax Theorem holds for $\Gi$ the proof concludes.

\section{Uniqueness: Proofs of Theorems~\ref{the3} and~\ref{the6}.}
\label{unique}
The proofs of Theorem~\ref{the3} and Theorem~\ref{the6}
will follow from three observations,
which we have written as separate lemmas. 
First, any limit point of $\{u_n\}$ solves (\ref{eq0.7}(i)) and
belongs to $\D^+_{\rho}$: this is shown in Lemmas~\ref{lem6}
and~\ref{lem6.0}.
Second, solutions of (\ref{eq0.7}(i)) in $\D^+_{\rho}$ satisfy
$\Gi(u,\mu)+\lr=0$ for any $\mu\in \M_{\rho}$: this
is shown in Lemma~\ref{lem6.1}. Third,
by convexity of $h\mapsto\Gi(\exp(h),\mu)$ shown in Proposition~\ref{prop4.1},
there is a unique solution of $\Gi(u,\mu)+\lr=0$ for any $\mu\in \M_{\rho}$:
this is shown in Lemma~\ref{lem6.2}.

\bl{lem6.0}
If $u\in\D_{\rho}$, $\int\!\! u d\nur=1$, and $u$ satisfies (\ref{eq0.7}(i)),
then $u$ is positive $\nur$-a.s.\ on $\A^c$.
\el
\bpr
We denote by $\B:=\{\eta:\ u(\eta)=0\}$. Since $u\in\D_{\rho}$,
we have for $i\not\in \S$ and $\eta$ $\nur$-a.s.,
\[
u(\eta)\ge u(\An_i^+ \eta)\quad\text{and}\quad
u(\An_i^+ \eta)\ge \frac{1}{1-\e_i} u(\eta).
\]
Thus, for $i\not\in \S$, $\B=(\An_i^+)^{-1}(\B)$ $\nur$-a.s.\ .
For any cylinder $\theta$ with base in $\N^{\S}\bs \A$, we will consider
$\B_{\theta}:=\B\cap \theta$. If $T^{i,j}$ denotes the
exchange operator at site $i,j\in \Z^d$, then
\[
\B_{\theta}
\stackrel{\scriptstyle{\nur-\text{a.s}}}{=}
(\An_i^+)^{-1}(\B_{\theta}),\quad\forall i\not\in\S
\Longrightarrow
\B_{\theta}
\stackrel{\scriptstyle{\nur-\text{a.s}}}{=}
(T^{i,j})^{-1}(\B_{\theta}),\quad\forall i,j\not\in\S.
\]
Indeed, 
\[
\B_{\theta}
\stackrel{\scriptstyle{\nur-\text{a.s}}}{=}
\bigcup_{k,l\in \N} \B_{\theta}\cap\{\eta(i)=k,\eta(j)=l\},
\]
so that we can go from 
\[
\B_{\theta}\cap\{\eta(i)=k,\eta(j)=l\}\quad\text{to}\quad
\B_{\theta}\cap\{\eta(i)=l,\eta(j)=k\}=(T^{i,j})^{-1}(\B_{\theta}
\cap\{\eta(i)=k,\eta(j)=l\}
\]
by a finite succession of creation
and annihilation of particles. Now, by Hewitt-Savage 0-1 law
on the lattice $\Z^d\bs \S$, we conclude that $\nur(\B_{\theta})
\in \{0,1\}$. Assume that for some cylinder $\theta$, $\nur(\B_{\theta})=1$.
Since $u$ satisfies (\ref{eq0.7}(i)) and $1_{\theta}\in \Lb$, we have
\be{eq6.17}
\int u \LL^*(1_{\theta})d\nur=0
\Longrightarrow
\sum_{i,j\in \Z^d} p^*(i,j)\int g_i(\eta) u(\eta)
1_{\theta}(T^i_j\eta) d\nur=0.
\ee
Now,
\be{eq6.18}
(T^i_j)^{-1}(\theta)=
\left\lbrace\begin{array}{l}
T^j_i(\theta)\text{ if } \theta(j)>0,\quad\text{and}\quad \emptyset
\text{ if } \theta(j)=0\quad\text{when}\quad i,j\in \S\\
\An_i^+(\theta)\quad\text{when}\quad i\in \S,\ j\not\in \S\\
\An_j^-(\theta)\text{ if } \theta(j)>0,\quad\text{and}\quad \emptyset
\text{ if } \theta(j)=0\quad\text{when}\quad i\not\in \S,\ j\in \S\\
\theta \quad\text{when}\quad i,j\not\in \S
\end{array}\right.
\ee
Since the moves on the right hand side generates all cylinders with base
in $\N^{\S}\bs \A$, we obtain 
\be{eq6.19}
\forall \theta\in \N^{\S}\bs \A,\quad\int_{\theta}\!\!\!u d\nur=0,
\ee
which is absurd since $\int ud\nur=1$. Thus, $\nur(\B)=0$ and
the proof is concluded.
\epr
\bl{lem6.1}
If $u\in \D_{\rho}$ satisfies (\ref{eq0.7}(i)) then
$\Gi(u,\mu)=-\lr$, for any $\mu\in \M_{\rho}$.
\el
\bpr
Let $u$ satisfies (\ref{eq0.7}(i)). By Lemma~\ref{lem6.0}, $u\in\D^+_{\rho}$.
For any $\v_n\in \Lb$ with $\v_n$ $\FF_n$-measurable,
we write (\ref{eq0.7}(i)) as
\be{eq6.2}
\int \LL^*(\v_n) ud\nur+\lr \int\v_n u d\nur=0.
\ee
We make the standard integration by parts and use
cancellations as in (\ref{eq4.9}) to obtain
\ba{eq6.4}
\int \LL^*(\v_n) u d\nur&=&
\sum_{i,j\in \L_n}\!\! \gr p^*(i,j)\int
\Delta^j_i \v_n \An_i^+ u d\nur
-\sum_{ i\in \L_n}\sum_{j\not\in \L_n}\!\!
\gr p^*(i,j)\int\nabla_i^+\v_n \An_i^+ u d\nur\cr
&&\qquad+\sum_{i\not\in \L_n}\sum_{j\in \L_n}\!\!
\gr p^*(i,j)\int\nabla_j^+\v_n \An_i^+ u d\nur\cr
&=&\sum_{i,j\in \L_n} \gr p^*(i,j)\int
\Delta^j_i(\v_n)\nabla_i^+(u)d\nur+\tilde R_n(\v_n),
\ea
where
\be{eq6.5}
\tilde R_n(\v_n)
=-\sum_{i\in \L_n}\gr p_n^*(i,i)\int \nabla_i^+\v_n \nabla_i^+u 
d\nur+ \sum_{i\not\in \L_n}\sum_{j\in \L_n}\gr p_n(j,i)
\int \nabla_j^+\v_n \nabla_i^+u d\nur.
\ee
Now, for any $\mu\in \M_{\rho}$ with density $f$, it is easy to note
that for a fix large integer $M$, 
\[
\v_n^{(M)}:=E_{\nur}[\frac{f}{u}\wedge M|\FF_n]\in \Lb,
\]
and if we set $\v=f/u$ and $\v^{(M)}=(f/u)\wedge M$, then both $\v^{(M)}$ and $\v$
are in $L^p(\nur)$ for any integer $p$, and are 
such that for $i$ large enough $|\nabla_i^+(\psi)|\le
2\psi(\e_i+\e_i^*)$. Indeed, for $i\not\in \S$ 
\be{eq6.8}
u\ge \An_i^+ u\ge u(1-\e_i),\quad\text{and}\quad
f\ge \An_i^+ f\ge f(1-\e_i-\e_i^*).
\ee
Thus, if $i$ is such that $1-\e_i-\e_i^*>0$,
\be{eq6.9}
\frac{f}{u}(1-\e_i-\e_i^*-1)\le \nabla_i^+(\frac{f}{u})\le
\frac{f}{u}(\frac{1}{1-\e_i}-1).
\ee
Thus, for $i$ large enough $|\nabla_i^+(\v)\le
2\v(\e_i+\e_i^*)$. Also, since $f\in L^p(\nur)$ and
$1_{\A^c}/u\in L^p(\nur)$ for any integer $p$ by Lemma~\ref{lem5.1},
we obtain that $\v\in L^p(\nur)$ for any $p$. The same is
true for $\v^{(M)}$ after a simple algebra.

By a reasoning by now standard, since $\v_n^{(M)}$ satisfies a
bound like (\ref{eq6.9})
\be{eq6.6}
|\tilde R_n(\v_n^{(M)})|\le c_1\sum_{i\in \p^R \L_n} (\e_i^2+(\e_i^*)^2)\int 
\v_n^{(M)} u d\nur
\le c_1 ||\v||_{\nur} ||u||_{\nur}
\sum_{i\in \p^R \L_n} \e_i^2+(\e_i^*)^2
\stackrel{\scriptstyle{n \ra \infty}}{\longrightarrow} 0.
\ee
Recall that for $i\in \L_n$, $\An_i^+\v_n^{(M)}=E_{\nur}[\An_i^
+\v^{(M)}|\FF_{\L_n}]$. Now, since $\v^{(M)}\in L^2(\nur)$, 
and $\{\v_n^{(M)},\ n\in \N\}$ is a positive martingale, we have
that $\{\v_n^{(M)}\}$ converges to $\v^{(M)}$ in $L^2(\nur)$ and a.s.\ .
Thus, for any $\psi\in L^2(\nur)$, and $i,j\in \L_n$
\be{eq6.3}
\lim_{n\to\infty} \int \v_n^{(M)} \psi d\nur=
\int \v^{(M)} \psi d\nur,\quad\text{and}\quad
\lim_{n\to\infty} \int \nabla_j^+\v_n^{(M)} \nabla_i^+\psi d\nur=
\int \nabla_j^+\v^{(M)} \nabla_i^+\psi d\nur.
\ee
Thus, combining (\ref{eq6.4}), (\ref{eq6.6}) and (\ref{eq6.3}) we obtain
(the series being absolutely convergent)
\be{eq6.7}
\sum_{i,j\in \Z^d}\gr p^*(i,j)\int \Delta^j_i(
\frac{f}{u}\wedge M)\nabla_i^+ u d\nur
+\lr \int (\frac{f}{u}\wedge M) u d\nur=0.
\ee
An identical expression to (\ref{eq6.7}) is
also valid for $f/u$ as we take the limit $M$ to infinity.

We will now show that $\Gi(u,\mu)$ has the same expression
as the first term of (\ref{eq6.7}).
Now, by taking the limit $n$ to infinity in expression (\ref{eq4.7bis}),
we obtain
\be{eq6.10}
\lim_{n\to\infty} \int \LL_n(u) \frac{f}{u} d\nur=
\sum_{i,j\in \Z^d}\gr p^*(i,j)\int \Delta^j_i \frac{f}{u} \nabla_i^+ u d\nur.
\ee
Indeed, (\ref{eq4.7bis}) only requires that $f/u\in L^p(\nur)$ and that
for $i$ large enough $|\nabla_i^+(\frac{f}{u})|\le
2\frac{f}{u}(\e_i+\e_i^*)$.  
Finally, since $\Gi(u,\mu)=\lim_{n\to\infty} \Gn(u,\mu)$, (\ref{eq6.7}) 
concludes the proof.
\epr
\bl{lem6.2} If $u,\tilde u\in \D^+_{\rho}$, and for any $\mu\in \M_{\rho}$
$\Gi(u,\mu)=\Gi(\tilde u,\mu)=-\lr$, and
$\int ud\nur=\int \tilde u d\nur$, then $u= \tilde u$ $\nur$-a.s.\ .
\el
\bpr
We can define
\[
h:=\log(u),\quad\text{and}\quad \tilde h:=\log(\tilde u),
\qquad\text{with}\qquad h,\tilde h\in \E_{\rho}. 
\]
Now, for $\gamma\in ]0,1[$, 
we form $h_{\g}=\g h+(1-\g)\tilde h$, and by convexity of $\tilde \Gn$,
for any $\mu\in \M_{\rho}$, 
\be{eq6.20}
0\le a_n(\mu):=\g {\tilde\Gn}(h,\mu)+(1-\g){\tilde\Gn}(\tilde h,\mu)-
{\tilde\Gn}(h_{\g},\mu)
\stackrel{\scriptstyle{n \ra \infty}}{\longrightarrow}
-\lr-{\tilde\Gi}(h_{\g},\mu),
\ee
where we used Lemma~\ref{lem6.1}. Now, Lemma~\ref{lem6.1}
is also valid for any $u^*$ limit point of $u^*_n$, the principal
eigenfunction of $\LL_n^*$. Note that since $u,\tilde u\in \D^+_{\rho}$
then $u,\tilde u\in L^2(\nur)$. By Jensen, this implies that
$\exp(h_{\gamma})\in L^2(\nur)$ and $\int \exp(h_{\gamma})u^*d\nur<
\infty$. Finally, Lemma~\ref{lem6.0} would
imply that $u^*|_{\A^c}>0$ $\nur$-a.s.\ , so that
$\int u^* \exp(h_{\g})d\nur>0$, and we can define
\be{eq6.22}
d\mu^*:=\frac{e^{h_{\g}}u^*d\nur}{\int\!\! e^{h_{\g}}u^*d\nur}\in\M_{\rho}.
\ee
Now, by duality, and Lemma~\ref{lem6.1} applied to $\LL^*$. 
\be{eq6.23}
\tilde \Gn(h_{\g},\mu^*)=\int \frac{\LL_n^*(u^*)}{u^*}  d\mu^*
=\Gamma_{\!\!n}^*(u^*,\mu^*)
\stackrel{\scriptstyle{n \ra \infty}}{\longrightarrow}
-\lr=\tilde \Gi(h_{\g},\mu^*).
\ee
Thus, $a_n(\mu^*)$ vanishes as $n$ tends to $\infty$. However, for any
$i,j\in \Z^d$ and $n$ large enough
\be{eq6.24}
a_n(\mu^*)\ge p(i,j)\int g_i A_{i,j} d\mu^*,\quad\text{with}\quad
0\le A_{i,j}:= \g e^{\nabla^i_j h}+(1-\g)e^{\nabla^i_j \tilde h}-
e^{(\g \nabla^i_j h+(1-\g)\nabla^i_j \tilde h)}.
\ee
Now $a_n(\mu^*)\to\infty$, $e^{h_{\g}}u^*>0$ $\nur$-a.s.\ on $\A^c$, and
(\ref{eq6.24}) imply that $p(i,j)g_i A_{i,j}=0$ $\nur$-a.s.\ on $\A^c$.
This in turn, implies that for $\eta(i)p(i,j)>0$, $\nur$-a.s., we have 
$\nabla^i_j h=\nabla^i_j \tilde h$ in $\A^c$.
Let us denote $f:=\tilde u/u$ on $\A^c$. Since, $p(.,.)$ is
irreducible, we obtain 
\be{eq6.25}
\forall i,j\quad\text{with }\eta(i)p(i,j)>0\qquad
f(T^i_j\eta)= f(\eta), 
\quad\nur-\text{a.s.\ }.
\ee
This in turn, implies that for $i,j\not\in\S$
$f(T^{i,j}\eta)=f(\eta)$ $\nur$-a.s.\ , so that
by Hewitt-Savage 0-1 law for exchangeable events, we conclude
that $f$ is $\nur$-a.s.\ constant on each cylinder $\theta$ with base in
$\N^{\S}\bs \A$, say $c_{\theta}:=f|_{\theta}$. 

We now show that the constants $\{c_{\theta}\}$ are the same.
Assume $\theta,\theta'\in \N^{\S}\bs \A$ with $T^i_j\theta=\theta'$.
If we denote $X_{\theta}:=f^{-1}(\{c_{\theta}\})$, then
\be{eq6.26}
\theta\subset X_{\theta},\qquad \theta'=T^i_j\theta\subset
(T^j_i)^{-1}(X_{\theta}),\qquad\text{and by (\ref{eq6.25})}\quad
(T^j_i)^{-1}(X_{\theta})
\stackrel{\scriptstyle{\nur\text{a.s}}}{=} X_{\theta}.
\ee
This yields $c_{\theta}=c_{\theta'}$. Assume now that for $j\in \S$ with
\[
\sum_{i\not\in \S}p(i,j)>0,\quad\text{we have}\quad
\theta'=\An_j^+ \theta.
\]
Take $i\not\in \S$ with $p(i,j)>0$, and note that
\be{eq6.27}
\theta'=\An_j^+ \theta\subset (T^j_i)^{-1}(X_{\theta})
\stackrel{\scriptstyle{\nur\text{a.s}}}{=} X_{\theta}
\quad[\text{by (\ref{eq6.25})}]
\ee
Thus, $c_{\theta}=c_{\theta'}$ in this case also. Now, we have
assumed that $\N^{\S}\bs \A$ was a connected set containing $0_{\S}$.
Thus, by a succession of moves $T^j_i$ and $\An_j^+$ applied to 
$0_{\S}$, we cover all of $\N^{\S}\bs \A$, and conclude that $f$ is
constant $\nur$-a.s.\ .
\epr
\section{Hitting time: Proof of Theorem~\ref{the4}.}
\label{hitting}
Let $u$ (resp. $u^*$) be the principal Dirichlet eigenfunction of
$\LL$ (resp. $\LL^*$) in $\D_{\rho}$ (resp. $\D_{\rho}^*$).
By Lemma~\ref{lem6.0}, $u$ and $u^*$ are $\nur$-a.s.\
positive on $\A^c$. Thus, we define a Markov semi-group on $\A^c$, 
\be{eq7.3}
\forall \eta\in \A^c,\qquad S_t^u(\v)(\eta):=
e^{\lr t}\frac{\bar S_t(u \v)}{u(\eta)}.
\ee
This semi-group is stationary with respect to 
\be{eq7.4}
d\hat{\mu}_{\rho}=\frac{u u^*d\nur}{\int\!\! u u^*d\nur}.
\ee
Note that since $1_{\A^c}/u\in L^2(\hat{\mu}_{\rho})$ by Lemma~\ref{app-key},
we have by definition, for all $\eta\in \A^c$,
\be{eq7.5}
S_t^u(\frac{1_{\A^c}}{u})(\eta)=
c_t\frac{u_t(\eta)}{u(\eta)}\quad\text{with}\quad
c_t=e^{\lr t}P_{\nur}(\tau>t),\quad\text{and}\quad
u_t:=\frac{P_{\eta}(\tau>t)}{P_{\nur}(\tau>t)}.
\ee
From (\ref{eq0.16}), $c_t\in [\underline{c},1]$, whereas
$u_t\in \D_{\rho}$ from inequality (4.7) of~\cite{ac}. 
It is then easy to check directly that, for any $t>0$
\be{eq7.6}
u S_t^u(\frac{1_{\A^c}}{u})\in \D_{\rho}\Longrightarrow
\text{if}\quad 
\psi_t:=
\frac{1}{t}\int_0^t\!\! S_s^u(\frac{1_{\A^c}}{u})ds,
\quad\text{then}\quad u\psi_t\in \D_{\rho}.
\ee
Now, by Jensen's inequality, $\{ S_t^u,t>0\}$ is a contraction
semi-group on $L^2(\A^c, \hat{\mu}_{\rho})$. Thus, by 
von Neumann's mean ergodic theorem in Hilbert space
(see e.g.\ \cite{petersen} Th.1.2 page 24), we obtain
\be{eq7.7}
\psi_t
\stackrel{\scriptscriptstyle{t\to\infty}}{\longrightarrow}
\psi\quad\text{in}\quad L^2(\hat{\mu}_{\rho}),
\qquad\text{and for any $t\ge 0$}\quad
S_t^u(\psi)=\psi,\ \hat{\mu}_{\rho}-\text{a.s.} .
\ee
If $\{\psi_t\}$ converge
to $\psi$ in $L^2(\hat{\mu}_{\rho})$, then $\{\psi_t u\}$ converge weakly
towards $\psi u$. Indeed, for any $\v$ bounded and continuous
\ba{eq7.8}
|\int \v u(\psi_t-\psi)d\nur|&\le&\left(\int\!\! u u^*d\nur\right)
||\frac{\v}{u^*}||_{\hat{\mu}_{\rho}}||\psi_t-\psi||_{\hat{\mu}_{\rho}}\cr
&\le& |\v|_{\infty} \left(\int\!\! u u^*d\nur\int\!\!\frac{u}{u^*}d\nur
\right)^{1/2}||\psi_t-\psi||_{\hat{\mu}_{\rho}}.
\ea
Since $\int \psi_t ud\nur=\int_0^tc_sds/t\le 1$,
the Remark~\ref{rem-comp} yields that $u\psi\in \D_{\rho}$.
Finally, since $\hat{\mu}_{\rho}$ and $\nur$ are equivalent in $\A^c$,
(\ref{eq7.7}) implies that for any $t\ge 0$,
\be{eq7.9}
\bar S_t(\psi u)=e^{-\lr t} \psi u\ \text{on }\A^c\ \nur-\text{a.s.}.
\ee
Thus, by differentiating (\ref{eq7.9}) at $t=0$, we obtain that
$\psi u$ is a Dirichlet principal eigenfunction in $\D^+_{\rho}$,
with $\underline{c} \le \int\!\! u\psi d\nur\le 1$.
By Theorem~\ref{the3}, this means that $\psi$ is constant. To
find the value of $\psi$, integrate (\ref{eq7.7}) against $1_{\A^c}$.
\be{eq7.10}
\psi\equiv \int \psi d\hat{\mu}_{\rho}=\lim_{t\to\infty}
\int\frac{1}{t}\int_{0}^t S_s^u(\frac{1}{u})ds\ d\hat{\mu}_{\rho}=
\int\frac{1}{u} d\hat{\mu}_{\rho}=\frac{1}{\int\!\! u u^* d\nur}.
\ee
Finally, since $1/u^*\in L^2(\hat{\mu}_{\rho})$, 
we integrate (\ref{eq7.7}) against $1/u^*$ to conclude the proof with
\be{eq7.11}
\frac{1}{t}\int_{0}^t\frac{e^{\lr s}P_{\nur}(\tau>s)}{
\int\!\! u u^* d\nur}ds=
\int\!\frac{1}{u^*}\psi_td\hat{\mu}_{\rho}
\stackrel{\scriptstyle{t \ra \infty}}{\longrightarrow}
\int\!\frac{1}{u^*}\psi d\hat{\mu}_{\rho}=\frac{1}{(\int\!\! u u^* d\nur)^2}.
\ee
\qed

\section{Appendix}
\label{stoch-domination}
We have often used Lemma~\ref{lem5.1} below to obtain regularity
of probability densities satisfying a gradient bound (\ref{def-D})
~\cite{ad1,af,ac,ad2}.
For ease of reading, we recall its simple proof. Then, 
in Lemma~\ref{app-key},
we show how similar arguments yield the regularity of $1_{\A^c}/\v$ for 
$\v\in\D^+_{\rho}$.

For ease of writing, we identify a cylinder with its base. 
Thus, when we write $\theta\in \N^{\S}$, we mean $\theta:=\{\eta\in
\N^{\Z^d}:\ \eta(i)=\theta(i),\ \forall i\in \S\}$.
Recalling the notations
used in the definition of $\nur$ (see (\ref{def-m})), let $\nu_{\epsilon}$ 
be the product measure
\[
d \nu_{\epsilon}(\eta) = \prod_{i \in \S} d\theta_{\gamma(\rho)}(\eta_i)
\prod_{i \notin \S} d\theta_{(1-\epsilon_i)\gamma(\rho) }(\eta_i)\, .
\]
We showed in \cite{ac} that when $d\ge 3$, $\nu_{\e}$ is abolutely
continuous with respect to $\nur$, and that if $\psi_{\e}:=
d\nu_{\e}/d\nur$, then for any integer $p$
\be{recall-psi}
\int \psi^p_{\e}d\nur<\infty,\qquad\text{and}\qquad
\int \frac{1}{\psi^p_{\e}}d\nur<\infty.
\ee
\br{rem-recall}
The purpose of introducing $\psi_{\e}$ was that for any
$i\not\in \S$, $\An_i^+\psi_{\e}=(1-\e_i)\psi_{\e}$. Thus, if
$\v\in \D_{\rho}$, then $\v/\psi_{\e}$ is increasing outside $\S$. Indeed,
using (\ref{def-Dn})(ii),
\[
\forall i\not\in \S,\qquad \An_i^+(\v/\psi_{\e})\ge \v/\psi_{\e}.
\]
\er

\bl{lem5.1} We assume that $d\ge 3$.  For any integer $n$, 
any $\theta\in \N^{\S}\bs\A$, and $\v\in \D_{\rho}$
\be{eq5.0}
\int_{\theta} \v^n d\nur\le \left( \frac{\int_{\theta}\v d\nur}
{\nu_{\e}(\theta)}\right)^n\int_{\theta} \psi_{\e}^n d\nur.
\ee
Also, 
\be{eq5.1}
\int\!\! \v^n d\nur\le C_n\left(\int\!\! \v d\nur\right)^n
\qquad\text{with}\qquad
C_n:=\frac{\int\!\! \psi_{\e}^n d\nur} {\nu_{\e}(0_{\S})^{n+1}}<\infty.
\ee
\el
\bpr
We define the measure $d\mu=\v d\nur$, and for $\theta\in \N^{\S}\bs\A$,
we define two probability measures $d\mu_{\theta}=1_{\theta}d\mu/\mu(\theta)$
and $d\nu_{\theta}=1_{\theta}d\nu_{\e}/\nu_{\e}(\theta)$.
Note that on $\theta$, the probability measure $\nu_{\e}$
satisfies Holley's condition (see Theorem 2.9, p.75 in~\cite{liggett})
which implies that it satisfies FKG's inequality.

\noindent{\bf Step 1.} We first show that for any $\phi$ 
decreasing on $\theta$,
\be{eq5.2}
\int \phi d\mu_{\theta} \le \int \phi d\nu_{\theta}.
\ee
By the Remark~\ref{rem-recall},
$d\mu_{\theta}/d\nu_{\theta}$ is increasing in $\theta$.
We apply FKG's inequality on $\theta$ 
\be{eq5.4}
\int \phi d\mu_{\theta}=\int \phi \frac{d\mu_{\theta}}{d\nu_{\theta}}
d\nu_{\theta}\le \int \phi d\nu_{\theta}.
\ee

\noindent{\bf Step 2.} First, note that $\v$ and
$\psi_{\e}=d\nu_{\e}/d\nur$ are non-negative decreasing on $\theta$.
So is $\v^i\psi_{\e}^j$ for any integers $i,j$. We apply (\ref{eq5.2})
to $\phi:=\v^i\psi_{\e}^j$ and obtain
\be{eq5.5}
\int_{\theta} \v^{i+1} \psi_{\e}^j \frac{d\nur}{\mu(\theta)}=
\int \v^i\psi_{\e}^j d\mu_{\theta} \le \int \v^i\psi_{\e}^j d\nu_{\theta}=
\int_{\theta} \v^{i} \psi_{\e}^{j+1} \frac{d\nur}{\nu_{\e}(\theta)}.
\ee
By induction, we obtain (\ref{eq5.0}) for any integer $n$.
Now, (\ref{eq5.1}) obtains after taking $\theta=0_{\S}$ and using
FKG's inequality once more. Indeed, since $\v$ and $1_{0_{\S}}$ are
both decreasing
\be{eq5.6}
\int \v^n d\nur\le \frac{\int \v^n 1_{0_{\S}} d\nur}{\nur(0_{\S})}.
\ee
\epr
\br{rem-M} Actually if $\mu\in\M_{\rho}$, then its density
$f:=d\mu/d\nur$ satisfies an inequality like (\ref{eq5.1})
but with $U:=\{i: (\e_i+\e^*_i)\ge 1\}$ replacing $\S$ which
was the domain where $\e_i=1$. Since $U$ is bounded,
$\nu_{\e}(U)>0$.
\er
\bl{app-key} We assume that $d\ge 3$.
Let $\v\in \D_{\rho}^+$ and $\theta\in \N^{\S}\bs\A$. Then,
for any integer $n$
\be{eq5.7}
\int_{\theta} \v^n d\nur\int_{\theta}\frac{1}{ \v^n} d\nur\le
\int_{\theta} \psi_{\e}^n d\nur \int_{\theta} \frac{1}{\psi_{\e}^n} d\nur.
\ee
Furthermore, for
\[
c_{\v,n}:=\sup_{\theta\in \N^{\S}\bs\A}\{\nur(\theta)
\left( \int_{\theta} \v \frac{d\nur}{\nur(\theta)}\right)^{-n}
\}<\infty,
\]
we have 
\be{eq5.8}
\int_{\A^c}\!\!\frac{1}{ \v^n} d\nur\le c_{\v,n}
\int\!\! \psi_{\e}^n d\nur \int_{\A^c}\!\! \frac{1}{\psi_{\e}^n} d\nur.
\ee
\el
\bpr
Recall that $\v/\psi_{\e}$ is increasing on $\A^c$
whereas for any integer $n$, $\v^n$ is decreasing. 
Thus, for any cylinder $\theta\in \N^{\S}\bs\A$, if we denote 
$d\tilde{\nur}=1_{\theta}d\nur/\nur(\theta)$, then, by FKG's inequality
\be{eq5.9}
\int_{\theta}\frac{1}{ \psi_{\e}^n} d\tilde \nur=
\int_{\theta}\left(\frac{\v}{\psi_{\e}}\right)^n\frac{1}{ \v^n} d\tilde \nur\ge
\int_{\theta}\left(\frac{\v}{\psi_{\e}}\right)^n d\tilde \nur
\int_{\theta}\frac{1}{ \v^n} d\tilde \nur.
\ee
Also, since $\psi_{\e}^n$ is decreasing
\be{eq5.10}
\int_{\theta}\left(\frac{\v}{\psi_{\e}}\right)^n d\tilde \nur
\int_{\theta} \psi_{\e}^n d\tilde \nur\ge
\int_{\theta}\left(\frac{\v}{\psi_{\e}}\right)^n \psi_{\e}^nd\tilde \nur
=\int_{\theta} \v^n d\tilde \nur
\ee
Multiplying (\ref{eq5.9}) by
$\int_{\theta} \psi_{\e}^n d\tilde \nur$, using (\ref{eq5.10}), and simplifying
by $\int_{\theta} (\v/\psi_{\e})^n d\tilde \nur>0$, (since $\v>0$ $\nur$-a.s.)
we obtain (\ref{eq5.7}).
Note that $c_{\v}<\infty$ since there is a finite number of elements in
$\N^{\S}\bs\A$, on each of which $\int_{\theta} \v d\nur>0$.

Finally, (\ref{eq5.8}) is obtained by summing over all $\theta\in
\N^{\S}\bs \A$, and applying H\"older's inequality to 
$\int\!\!\v^n d\tilde{\nur}$.
\epr

\noindent {\bf Acknowledgements}.We thank 
Fabienne Castell for her comments on a first version of this paper.

\end{document}